\documentclass[11pt]{amsart}

\usepackage[all]{xy}
\usepackage{color, pstricks}
\usepackage{mathrsfs}
\usepackage{amsmath,amsthm, amssymb, amsgen,amsxtra,amsfonts, amsbsy} 
\usepackage[all]{xy}
\usepackage{verbatim}
\usepackage{stmaryrd}
\SetSymbolFont{stmry}{bold}{U}{stmry}{m}{n}
\usepackage{bm}
\usepackage{tikz,caption,float}
\usepackage[backref=page]{hyperref}
\usepackage{listings}
\usetikzlibrary{decorations.pathmorphing}
\usetikzlibrary{arrows,automata,positioning}

\hypersetup{
     colorlinks   = true,
     citecolor    = red,
     linkcolor = blue,
     urlcolor = purple
}
\theoremstyle{definition}
\newtheorem{Definition}{Definition}[section]

\newtheorem{Remark}[Definition]{Remark}

\newtheorem*{Acknowledgements}{Acknowledgements}
\newtheorem*{AssumptionStar}{Assumption $(\star)$}

\theoremstyle{plain}
\newtheorem{Theorem}[Definition]{Theorem}

\newtheorem{Proposition}[Definition]{Proposition}

\newtheorem{Lemma}[Definition]{Lemma}
\newtheorem{Corollary}[Definition]{Corollary}

\newtheoremstyle{voiditstyle}{3pt}{3pt}{\itshape}{\parindent}%
{\bfseries}{.}{ }{\thmnote{#3}}%
\theoremstyle{voiditstyle}

\newtheoremstyle{voidromstyle}{3pt}{3pt}{\rm}{\parindent}%
{\bfseries}{.}{ }{\thmnote{#3}}%
\theoremstyle{voidromstyle}

\newcommand{\prf}[1][Proof]{\par\noindent{\sc #1.}\quad}

\newcommand{\cal}{\mathcal}

\newcommand{\DD}{{\mathbb{D}}}

\newcommand{\ZZ}{{\mathbb{Z}}}

\newcommand{\QQ}{{\mathbb{Q}}}

\newcommand{\iso}{\cong}

\newcommand{\OO}{{\cal O}}

\newcommand{\onto}{{\twoheadrightarrow}}
\newcommand{\Spec}{{\rm Spec}\:}

\newcommand{\BK}{{\rm BK}}

\newcommand{\Bcris}{B_{{\rm cris}}}
\newcommand{\Bst}{B_{{\rm st}}}
\newcommand{\BdR}{B_{{\rm dR}}}
\newcommand{\DDcris}{\DD_{{\rm cris}}}
\newcommand{\DDst}{\DD_{{\rm st}}}
\newcommand{\DDdR}{\DD_{{\rm dR}}}
\newcommand{\DDpst}{\DD_{{\rm pst}}}

\newcommand{\Het}[1]{{H_{\rm \acute{e}t}^{{#1}}}}

\newcommand{\Hcris}[1]{{H_{\rm cris}^{{#1}}}}

\newcommand{\HdR}[1]{{H_{\rm dR}^{{#1}}}}

\newcommand{\et}{\mathrm{\acute{e}t}}
\newcommand{\pow}[1]{\llbracket #1 \rrbracket}	
\newcommand{\isomto}{\overset{\sim}{\rightarrow}}

\title[A N\'eron--Ogg--Shafarevich Criterion for K3 Surfaces]{A N\'eron--Ogg--Shafarevich Criterion for K3 Surfaces\footnote{\lowercase{2010 \uppercase{M}athematics \uppercase{S}ubject \uppercase{C}lassification: 14\uppercase{J}28 (primary), 11\uppercase{G}25, 14\uppercase{F}20, 14\uppercase{F}30, 14\uppercase{G}20 (secondary)}}}

\author{Bruno Chiarellotto}
\address[Chiarellotto]{Universit\`a degli Studi di Padova, Dipartimento di Matematica ``Tullio Levi-Civita'', Via Trieste 63, 35121 Padova, Italy}
\email{chiarbru@math.unipd.it}

\author{Christopher Lazda}
\address[Lazda]{Universiteit van Amsterdam, Korteweg--de\thinspace Vries Institute for Mathematics, P.O. Box 94248, Amsterdam 1090 GE, Netherlands}
\email{c.d.lazda@uva.nl}

\author{Christian Liedtke}
\address[Liedtke]{TU M\"unchen, Zentrum Mathematik - M11, Boltzmannstr. 3, 85748 Garching bei M\"unchen, Germany}
\email{liedtke@ma.tum.de}

\begin{document}

\begin{abstract}
The naive analogue of the N\'eron--Ogg--Shafarevich criterion is false for K3 surfaces, that is, there exist K3 surfaces over Henselian, discretely valued fields $K$, with unramified $\ell$-adic \'etale cohomology groups, but which do not admit good reduction over $K$. Assuming potential semi-stable reduction, we show how to correct this by proving that a K3 surface has good reduction if and only if $\Het{2}(X_{\overline{K}},\QQ_\ell)$ is 
unramified, and the associated Galois representation over the residue field coincides with the second cohomology of a certain ``canonical reduction'' of $X$. We also prove the corresponding results for $p$-adic \'etale cohomology.
 \end{abstract}

\maketitle

\setcounter{tocdepth}{1}
\tableofcontents

\section{Introduction}

Let $\OO_K$ be a Henselian DVR,
whose residue field $k$ is perfect of characteristic $p\geq 0$, and whose fraction field $K$ is of characteristic $0$. 
Let $G_K={\rm Gal}(\overline{K}/K)$ be the absolute Galois group of $K$.
If $p>0$, let $W=W(k)$ be the Witt ring and $K_0=\mathrm{Frac}(W)$ its fraction field, thus we have $K_0\subset \widehat{K}$.

\subsection{Good reduction and Galois-representations}
Given a variety $X$ that is smooth and proper over $K$, one can ask whether
$X$ has {\em good reduction}, that is, whether there exists a smooth and proper 
model of $X$ over $\Spec\OO_K$.
If it exists, then
 \begin{enumerate} 
 \item the $G_K$-representation on $\Het{n}(X_{\overline{K}},\QQ_\ell)$
 is unramified for all $n$ and all primes $\ell\neq p$, and
 \item if $p>0$, the $G_K$-representation $\Het{n}(X_{\overline{K}},\QQ_p)$
 is crystalline for all $n$.
 \end{enumerate} 
This yields necessary conditions for good reduction of $X$ in terms of
$G_K$-representations, and we refer to Section \ref{sec: generalities}
for details.

\subsection{Abelian varieties}
If $A$ is an Abelian variety over $K$, then a famous theorem of Serre and Tate \cite{ST68}, 
extending previous work of N\'eron, Ogg, and Shafarevich on elliptic curves, states that
$X$ has good reduction 
if and only if the $G_K$-representation on $\Het{1}(A_{\overline{K}},\QQ_\ell)$
is unramified for one (resp. all) primes $\ell\neq p$.
If $p>0$, then good reduction of $A$ is also equivalent to the $G_K$-representation 
on $\Het{1}(A_{\overline{K}},\QQ_p)$ being crystalline, which was proven 
in special cases by Fontaine \cite{Fon79} and 
Mokrane \cite{Mok93} and by Coleman and Iovita \cite{CI99} 
in general.

\subsection{Curves} 
If $C/K$ is a curve (smooth, projective and geometrically connected), then results of 
Oda \cite[Theorem 3.2]{Oda95} show that $C$ has good reduction if and only if the outer $G_K$-action 
on its $\QQ_\ell$-unipotent fundamental group $\pi_1^\et(C_{\overline{K}})_{\QQ_\ell}$ 
is unramified for one (resp. all) primes $\ell\neq p$. 
The $p$-adic version of this result was proven by Andreatta, Iovita, and Kim \cite{AIK15}. 

\subsection{Models of K3 surfaces}
An interesting, rich, and accessible class of algebraic varieties
beyond curves and Abelian varieties are K3 surfaces.
We refer to Section \ref{sec: crys rep} for definitions.
To establish criteria of good reduction, it is oftentimes helpful to have 
good models of these varieties at hand and sometimes, it even suffices
to have these models only after base change to a finite extension of $K$.
For example, for curves, these are semi-stable models and for
Abelian varieties, these are N\'eron models 
(maybe even with semi-Abelian reduction).

For K3 surfaces, such good models were introduced
by Kulikov \cite{Kul77} and later studied by 
Pinkham--Persson \cite{PP81} and Nakkajima \cite{Nak00}.
More precisely, if $X$ is a K3 surface over $K$, then
a \emph{Kulikov model} is a flat and proper algebraic space
$$
   {\cal X}\,\to\,\Spec\OO_K
$$
with generic fiber $X$, whose special fiber is a strict
normal crossing divisor, and whose relative 
canonical divisor is trivial.
Following the discussion in \cite[Section 3.1]{LM18},
we make the following assumption.

\begin{AssumptionStar}
  A K3 surface $X$ over $K$ satisfies ($\star$) if
  there exists a finite field extension $L/K$ such that
  $X_L$ admits a Kulikov model.
\end{AssumptionStar}

Roughly speaking, ($\star$) would follow from
potential semi-stable reduction of K3 surfaces, which is
not known yet.
At the moment, ($\star$) is known to hold
in equal characteristic zero or in mixed characteristic $(0,p)$
if the K3 surface admits an ample line bundle 
of self-intersection strictly less than $(p-4)$.
We also note that Kulikov models usually only exist
as algebraic spaces and that they are usually not unique - 
even in the case of good reduction.
We refer to \cite{LM18} for precise
statements, additional results, and discussion.

\subsection{A \texorpdfstring{$\bm{p}$}{p}-adic criterion for potential good reduction of K3 surfaces}
In this article, we study criteria for good reduction of a K3 surface
$X$ over $K$ in terms of the $G_K$-representation 
on $\Het{2}(X_{\overline{K}},\QQ_\ell)$.
If $\ell\neq p$, then such a criterion was already established in 
\cite{LM18}.
Our first result is the following extension to the $p$-adic case.

\begin{Theorem}
 \label{thm: first main result in introduction}
  Suppose that $p>0$, and let $X$ be a K3 surface over $K$ that satisfies $(\star)$.
  Then, the following are equivalent:
   \begin{enumerate} 
    \item $X$ has good reduction after a finite and unramified extension of $K$. \label{num :1st main in intro 1}
    \item The $G_K$-representation on  $\Het{2}(X_{\overline{K}},\QQ_\ell)$ is
    unramified for one $\ell\neq p$.
    \item The $G_K$-representation on  $\Het{2}(X_{\overline{K}},\QQ_\ell)$ is
    unramified for all $\ell\neq p$.
    \item The $G_K$-representation on  $\Het{2}(X_{\overline{K}},\QQ_p)$ is
    crystalline.
   \end{enumerate} 
\end{Theorem}

\begin{Remark} Since there is a canonical isomorphism $G_K\cong G_{\widehat{K}}$ between the absolute Galois group of $K$ and that of its completion, it makes sense to speak of $\Het{2}(X_{\overline{K}},\QQ_p)$ being crystalline even if $K$ is not necessarily complete.
\end{Remark}

The equivalences $(1)\Leftrightarrow(2)\Leftrightarrow(3)$ were already
established in \cite{LM18} and Matsumoto \cite{Mat15} showed
that when $K$ is complete, $(4)$ implies good reduction after some finite, but possibly ramified,
extension of $K$.

\begin{Remark}
 \label{rem: counter examples}
 Using the examples constructed in \cite{LM18},
 we see that the implication $(4)\Rightarrow(1)$ requires in general
 a {\em non-trivial} unramified extension:
 more precisely, for every prime $p\geq5$,
 there exists a K3 surface over $K=\QQ_p$ that
 satisfies conditions (2), (3), and (4) of Theorem \ref{thm: first main result in introduction}
 but that does \emph{not} have good reduction over $K$, 
 see Theorem \ref{thm: counter examples}.
\end{Remark}

\subsection{RDP models and the canonical reduction}

Thus, although we usually do not have good reduction over $K$ in the situation of
Theorem \ref{thm: first main result in introduction}, we may appeal to \cite{LM18} 
to construct certain ``RDP models'' of K3 surfaces 
over $K$, which are models that have reasonably mild singularities.

More precisely, let $X$ be a K3 surface over $K$ that satisfies $(\star)$ and
the equivalent conditions of Theorem \ref{thm: first main result in introduction}.
Then, there exists a projective and flat model
$$ 
  \cal{X} \,\rightarrow\, \Spec \OO_K
$$
for $X$ such that the special fiber $\cal{X}_k$ has at worst canonical singularities, and
its minimal resolution is a K3 surface over $k$, see Theorem \ref{thm: RDP models}.
We refer to Definition \ref{def: singular surfaces} for the notion of canonical singularities. 

It follows from a theorem of Matsusaka and Mumford \cite{MM64} that the minimal resolution of singularities 
$$
Y \,\to\, \cal{X}_k
$$
of the special fiber is unique up to unique isomorphism. In particular, it does not depend on the choice of RDP model $\cal{X}$.
If $X$ does in fact have good reduction over $K$, then $Y$ is simply the special fiber 
of any smooth model. We will refer to $Y$ as the \emph{canonical reduction} of $X$.

It is important to stress that the models $\cal{X}$ asserted by
Theorem \ref{thm: first main result in introduction} and
Theorem \ref{thm: RDP models} are usually not unique.
However, once a polarization, that is an ample line bundle $\cal{L}$ on $X$, is chosen,
there is a canonical choice of an RDP model
$$ 
 P(X,\cal{L}) \,\to\, \Spec \OO_K
$$
that we call the \emph{canonical RDP model} of the pair $(X,\cal{L})$,
see the discussion in Section \ref{sec: rdpcr}. 

 \subsection{The Weyl group of a polarized K3 surface}
 The existence of examples as in Remark \ref{rem: counter examples} raises the 
question of how to distinguish between good reduction, and good reduction 
over an unramified extension. 
To do so, we introduce the Weyl group of a polarized K3 surface $(X,\cal{L})$ 
over $K$ that has good reduction over an unramified extension.

Let $P(X,\cal{L})\rightarrow \Spec\OO_K$ be the canonical RDP model of the pair $(X,\cal{L})$.
Then, its special fiber has canonical singularities, and the minimal resolution of singularities
$$ 
   Y \,\to\, P(X,\cal{L})_k
$$
is given by the canonical reduction of $X$. Denote by $E_{X,\cal{L}}\subset Y$ the exceptional locus.

Over $\bar{k}$ the divisor $E_{X,\cal{L}}$ splits into a union of smooth and rational curves,
which have self-intersection number $-2$.
Then, the \emph{Weyl group} of $(X,\cal{L})$ is the subgroup
$$
 \cal{W}_{X,\cal{L}}^\mathrm{nr} \,\leq \, \mathrm{Aut}_{\ZZ}\left(\mathrm{Pic}(Y_{\bar k})\right)
$$
that is generated by the reflections in these $(-2)$-curves of $E_{X,\cal{L},\bar{k}}$, we refer to Section \ref{subsec: RDP and K3} for details.
The superscript ``$\mathrm{nr}$'' refers to the fact that we are working over $\bar{k}$, 
or, equivalently over $K^\mathrm{nr}$, and there are similar versions over any finite unramified extension $L/K$. 
We note that the absolute Galois group $G_k=\mathrm{Gal}(\bar{k}/k)$ naturally acts on 
$\cal{W}_{X,\cal{L}}^\mathrm{nr}$ via its action on $\mathrm{Pic}(Y_{\bar k})$, or, equivalently on $E_{X,\cal{L},\bar k}$. 

We now come to our first improvement of Theorem \ref{thm: first main result in introduction},
which yields a necessary and sufficient criterion for good reduction of a 
K3 surface over $K$.
This result also ``explains'' the counter-examples from \cite{LM18}.

\begin{Theorem} 
  \label{theo: intro gr} 
  Let $X$ be a K3 surface over $K$ that satisfies $(\star)$ and
  the equivalent conditions of Theorem \ref{thm: first main result in introduction}.
  Then, for any choice $\cal{L}$ of polarization on $X$,
  there exists a non-abelian cohomology class
  $$
    \alpha_{X,\cal{L}}^\mathrm{nr} \,\in\, H^1(G_k,\cal{W}_{X,\cal{L}}^\mathrm{nr})
 $$
  depending only on $(X,\cal{L})$, which is trivial if and only if 
  $X$ has good reduction over $K$.
 \end{Theorem}

\subsection{Another interpretation of the cohomology class} \label{intro:cohclass}
 As well as controlling the failure of $X$ to have good reduction over $K$ itself, 
the cohomology class $\alpha_{X,\cal{L}}^\mathrm{nr}$ also describes the difference between the $G_k$-modules
$$
       \Het{2}(X_{\overline{K}},\QQ_\ell),\mbox{ \quad and \quad }\Het{2}(Y_{\bar{k}},\QQ_\ell) 
$$
for all primes $\ell\neq p$, as well as the difference between the $F$-isocrystals
$$ 
     \DDcris\left(\Het{2}(X_{\overline{K}},\QQ_p)\right)\mbox{ \quad and \quad }\Hcris{2}(Y/K_0)
$$
when $p>0$, where $Y$ denotes the canonical reduction of $X$ (again, we are using the fact that $G_K\cong G_{\widehat{K}}$ to define $\DDcris$ in the latter case). More precisely, there are $G_k$-equivariant homomorphisms
$$
\begin{array}{lclll}
  \cal{W}_{X,\cal{L}}^\mathrm{nr} &\to& \mathrm{Aut}_{\QQ_\ell} & \left(\Het{2}(Y_{\bar{k}},\QQ_\ell) \right)  & (\ell\neq p)\\
  \cal{W}_{X,\cal{L}}^\mathrm{nr} &\to& \mathrm{Aut}_{K_0^\mathrm{nr},F}&\left(\Hcris{2}(Y/K_0) \otimes K_0^\mathrm{nr}\right) & (p>0)
\end{array}
$$
where $K_0^\mathrm{nr}$ is the maximal unramified extension of $K_0$, and the subscript ``$F$'' in the refers to the fact that we are considering Frobenius equivariant automorphisms. 
Using these homomorphisms and the induced map on cohomology, we obtain $\ell$-adic and crystalline 
realizations
$$
 \begin{array}{lcll}
   \beta_{X,\ell}^\mathrm{nr} &\in& 
    H^1\left(G_k,  \mathrm{Aut}_{\QQ_\ell}  ( \Het{2}(Y_{\bar{k}},\QQ_\ell) )  \right) & (\ell\neq p) \\
  \beta_{X,p}^\mathrm{nr} &\in& 
    H^1\left(G_k,   \mathrm{Aut}_{K_0^\mathrm{nr},F}  ( \Hcris{2}(Y/K_0) \otimes K_0^\mathrm{nr})  \right) & (p>0)
 \end{array}
$$
of the cohomology class $\alpha_{X,\cal{L}}^\mathrm{nr}$ from Theorem \ref{theo: intro gr} (we will see that these do not depend on $\mathcal{L}$). 

By the general theory of descent, we can then ``twist'' $\Het{2}(Y_{\bar{k}},\QQ_\ell)$ and $\Hcris{2}(Y/K_0)$ by these cohomology classes
to obtain new $G_k$-modules and, when $p>0$, $F$-isocrystals over $k$,
$$
\Het{2}(Y_{\bar{k}},\QQ_\ell)^{\beta^\mathrm{nr}_{X,\ell}} \mbox{ \quad and \quad }
 \Hcris{2}(Y/K_0)^{\beta^\mathrm{nr}_{X,p}},
$$
respectively. After these preparations, we can now state our next result.

\begin{Theorem} 
 \label{theo: intro descr} Under the assumptions and notations of  Theorem \ref{theo: intro gr},
 there are natural isomorphisms 
$$
\begin{array}{lcll}
    \Het{2}(X_{\overline{K}},\QQ_\ell) &\cong& \Het{2}(Y_{\overline{k}},\QQ_\ell)^{\beta^\mathrm{nr}_{X,\ell}} & (\ell\neq p) \\
    \DDcris(\Het{2}(X_{\overline{K}},\QQ_p)) &\cong& \Hcris{2}(Y/K_0)^{\beta^\mathrm{nr}_{X,p}} & (p>0)
  \end{array}
$$
of $G_k$-modules and $F$-isocrystals over $k$, respectively.
\end{Theorem}

\begin{Remark} Note that this result also implies that the classes $\beta_{X,\ell}^\mathrm{nr}$ and $\beta_{X,p}^\mathrm{nr}$ are independent of the choice of polarisation $\mathcal{L}$.
\end{Remark}

\subsection{A N\'eron--Ogg--Shafarevich criterion for K3 surfaces}

Our main result is then the following analogue of the N\'eron--Ogg-Shafarevich criterion for the good reduction of K3 surfaces over $K$.

\begin{Theorem} 
  \label{theo: intro NOS} 
  Let $X$ be a K3 surface over $K$ satisfying $(\star)$. 
  Then, the following are equivalent.
   \begin{enumerate} 
    \item $X$ has good reduction over $K$.
    \item There exists a prime $\ell\neq p$ such that (resp. for all primes $\ell\neq p$) 
     the $G_K$-representation on $\Het{2}(X_{\overline{K}},\QQ_\ell)$ is unramified 
     and there exists an isomorphism
     $$ 
        \Het{2}(X_{\overline{K}},\QQ_\ell) \,\cong\, \Het{2}(Y_{\bar{k}},\QQ_\ell)
     $$
     of $G_k$-modules.
    \item If $p>0$, $\Het{2}(X_{\overline{K}},\QQ_p)$ is crystalline and there exists an isomorphism
     $$
         \DDcris(\Het{2}(X_{\overline{K}},\QQ_p)) \,\cong\, \Hcris{2}(Y/K_0)
     $$
     of $F$-isocrystals over $k$.
 \end{enumerate} 
Here, $Y$ denotes the canonical reduction of $X$. \end{Theorem}

By Theorem \ref{theo: intro gr} and Theorem \ref{theo: intro descr}, the proof of this amounts to showing that a certain map in non-abelian group cohomology has trivial kernel. The details of this reduction are spelled out in Section \ref{sec: NOS}, and the group cohomology calculation is performed in Section \ref{sec: nagc}.

\subsection{Integral \texorpdfstring{$\bm{p}$}{p}-adic Hodge theory}
In fact, these results hold integrally as well, although the statement of this in the $p$-adic case needs some set-up.
Indeed, if $p>0$, if
$$
  \rho \,:\,G_K \,\to\, \mathrm{GL}(V)
$$
is a crystalline $G_K$-representation, and if $\Lambda\subseteq  V$ 
is a $G_K$-stable $\ZZ_p$-lattice, then Kisin \cite{Kis06} associated to it a
Breuil--Kisin module $\BK_{\OO_K}(\Lambda)$, which is a $\mathfrak{S}:=W\pow{u}$-module 
together with some extra data (we refer to
Section \ref{sec: generalities} for more details). 
In particular, we may specialize such an object along the homomorphism $\mathfrak{S}\rightarrow W$ given by the composite of the ``$u=0$'' map and the Frobenius on $W$, to obtain an $F$-crystal 
\[ \BK_{\OO_K}(\Lambda) \otimes_{\mathfrak{S},\sigma} W \]
over $W$ (the reason for using this map rather than simply the ``$u=0$'' map is explained in Section \ref{sec: generalities}).

If $X$ is a smooth and proper variety over $K$ with good reduction,
say via some smooth and proper model ${\cal X}\to\Spec\OO_K$
with special fiber ${\cal X}_k$.
Let us assume moreover that ${\cal X}_k$ is a scheme and that 
$\Hcris{\ast}({\cal X}_k/W)$ is torsion-free.
Then, the $G_K$-representations on 
$\Het{n}(X_{\overline{K}},\QQ_p)$ are crystalline,
and after results of Bhatt, Morrow, and Scholze  \cite{BMS15,BMS16}, 
there exists an isomorphism of $F$-crystals
$$
   \BK_{\OO_K}\left( \Het{n}(X_{\overline{K}},\ZZ_p) \right) \otimes_{\mathfrak{S},\sigma} W
     \,\cong\, \Hcris{n}({\cal X}_k/W).
$$
which is compatible with the usual (rational) crystalline 
comparison theorem. 
We also refer to Section \ref{sec: generalities} for details.

Now, suppose that $X$ is a K3 surface over any Henselian $K$ as above (i.e. allowing $p=0$), that satisfies $(\star)$ and
  the equivalent conditions of Theorem \ref{thm: first main result in introduction}.
Then we obtain the same relationship between 
$$\Het{2}(X_{\overline{K}},\ZZ_\ell) \;\;\text{and}\;\; \Het{2}(Y_{\bar k},\ZZ_\ell),$$ as well as between 
$$ \BK_{\OO_K}\left( \Het{2}(X_{\overline{K}},\ZZ_p) \right) \otimes_{\mathfrak{S},\sigma} W\;\; \text{and} \;\; \Hcris{2}(Y/W),$$
as we saw in Theorem \ref{theo: intro descr}  
in the rational case. 
Namely, there are $G_k$-equivariant homomorphisms
$$
\begin{array}{lclll}
  \cal{W}_{X,\cal{L}}^\mathrm{nr} &\to& \mathrm{Aut}_{\ZZ_\ell} & \left(\Het{2}(Y_{\bar{k}},\ZZ_\ell) \right)  & (\ell\neq p)\\
  \cal{W}_{X,\cal{L}}^\mathrm{nr} &\to& \mathrm{Aut}_{W^\mathrm{nr},F}&\left(\Hcris{2}(Y/W) \otimes W^\mathrm{nr}\right) & (p>0)
\end{array}
$$
where $W^\mathrm{nr}$ is the ring of integers in $K_0^\mathrm{nr}$. We therefore obtain integral realisations of the class $\alpha_{X,\cal{L}}^\mathrm{nr}$
$$ \begin{array}{lcll}
   \beta_{X,\ell}^\mathrm{nr} &\in& 
    H^1\left(G_k,  \mathrm{Aut}_{\ZZ_\ell}  ( \Het{2}(Y_{\bar{k}},\ZZ_\ell) )  \right) & (\ell\neq p)\\
  \beta_{X,p}^\mathrm{nr} &\in& 
    H^1\left(G_k,   \mathrm{Aut}_{K_0^\mathrm{nr},F}  ( \Hcris{2}(Y/W) \otimes W^\mathrm{nr})  \right) & (p>0)
 \end{array}
 $$
as before.

\begin{Theorem}
 \label{thm: last main result in introduction}
  Let $X$ be a K3 surface over $K$ that satisfies $(\star)$ and
   the equivalent conditions of Theorem \ref{thm: first main result in introduction}.
   \begin{enumerate} 
    \item  There are natural isomorphisms 
$$
 \begin{array}{lcll}
    \Het{2}(X_{\overline{K}},\ZZ_\ell) &\cong& \Het{2}(Y_{\overline{k}},\ZZ_\ell)^{\beta^\mathrm{nr}_{X,\ell}} & (\ell\neq p)\\
    \BK_{\OO_K}\left(\Het{2}(X_{\overline{K}},\ZZ_p)\right) \otimes_{\mathfrak{S},\sigma} W &\cong& \Hcris{2}(Y/W)^{\beta^\mathrm{nr}_{X,p}} & (p>0)
  \end{array}
$$
of $G_k$-modules and $F$-crystals over $k$, respectively.
    \item $X$ has good reduction over $K$ if and only if there exists an isomorphism
    $$ \Het{2}(X_{\overline{K}},\ZZ_\ell) \cong \Het{2}(Y_{\overline{k}},\ZZ_\ell) $$
    of $G_k$-modules for some $\ell\neq p$ (resp. for all $\ell\neq p$). When $p>0$, $X$ has good reduction over $K$ if and only if
     there exists an isomorphism
    $$
       \BK_{\OO_K}\left( \Het{2}(X_{\overline{K}},\ZZ_p) \right)  \otimes_{\mathfrak{S},\sigma} W
        \,\cong\,  \Hcris{2}(Y/W)
     $$
    of $F$-crystals over $k$.
   \end{enumerate} 
\end{Theorem}

\subsection{Equal and positive characteristic}
One might ask whether all these results also hold 
if $\OO_K$ is an excellent and Henselian DVR of equal characteristic $p>0$.
The main problem is Assumption $(\star)$: if we knew that
the extension $L/K$ there can always be chosen \emph{separable},
then the statements and proofs of \cite{LM18} and this article would most likely go 
through verbatim.
At the moment, however, we do not know whether it is reasonable
to expect this separability in general. We plan to return to this question in the future.

\begin{Acknowledgements}

We would like to thank two anonymous referees, the first of whom pointed out a serious error in a previous article of ours on crystalline Galois representations associated to K3 surfaces. This spurred us on to obtain (in our view) much more interesting results than we had originally envisaged, and this article would not exist without them. We would also like to greatly thank a second anonymous referee for a careful and detailed reading of the current article, and for the many helpful suggestions and comments  made.

 We also thank Yves Andr\'e, Fabrizio Andreatta, Bhargav Bhatt, Luc Illusie,
 Teruhisa Koshikawa, Remke Kloosterman, Martin Olsson and Peter Scholze for discussions and help.
 We especially thank Yuya Matsumoto for many comments on the first
 version of this article.
 The third named author would like to thank the two first named authors for an
 invitation to the University of Padova in October 2016.

The first named author is supported by the grant 
 MIUR-PRIN 2015 ``Number Theory and Arithmetic Geometry''. 
 The second named author was supported by a Marie Curie fellowship 
 of the Istituto Nazionale di Alta Matematica ``F. Severi'', and the Netherlands Organisation for Scientific Research
(NWO). The first and second named authors were supported by the Universit\`a di Padova grant PRAT 2015 ``Vanishing Cycles, Irregularity and Ramification''.
 The third named author is supported by the ERC Consolidator Grant
 681838 ``K3CRYSTAL''.
 \end{Acknowledgements}

\addtocontents{toc}{\setcounter{tocdepth}{0}}

\section*{Notations and Conventions}

\addtocontents{toc}{\setcounter{tocdepth}{1}}

Throughout the whole article, we fix the following notations
$$
\begin{array}{ll} 
  \OO_K & \mbox{a Henselian DVR} \\
  \pi & \mbox{a uniformizer of }\OO_K \\
  K & \mbox{the field of fractions of }\OO_K,\mbox{  which we assume to be}\\
  & \mbox{ \quad of characteristic }0 \\
  k & \mbox{the residue field of }\OO_K,\mbox{  which we assume to be}\\
  & \mbox{ \quad perfect of characteristic }p\geq 0\\
   \widehat{K} & \mbox{the completion of $K$} \\
  W=W(k) & \mbox{if $p>0$, the ring of Witt vectors of } k \\
  K_0 & \mbox{if $p>0$, the field of fractions of }W \\
  K_0^{{\rm nr}},\;K^{{\rm nr}} & \mbox{the maximal unramified extensions of }K_0\mbox{ and }K \\ & \quad \mbox{ respectively}\\
  \overline{K},\,\bar{k} & \mbox{a fixed algebraic closure of }K,\mbox{and its residue field, an} \\ & \mbox{\quad algebraic closure of }k \\
  G_K,\, G_k & \mbox{the absolute Galois groups of $K$ and $k$, respectively} 
 \end{array}
$$
If $L/K$ is a field extension, and $X$ is a scheme over $K$, we abbreviate
the base-change $X\times_{\Spec K}\Spec L$ by $X_L$. 
If $p>0$, and $Y$ is a scheme over $k$, we will write
$$
  \Hcris{n}(Y/K_0) \,:=\, \Hcris{n}(Y/W)\otimes_W K_0
$$
and similarly over any finite extension of $k$.

\section{Generalities}
\label{sec: generalities}

In this section, we recall a couple of general facts on models of varieties,
crystalline Galois representations, the functors of Fontaine and Kisin,
and their relation to good reduction.

\subsection{Models} 
\label{subsec: models}
We start with the definition of various types of models.

\begin{Definition}
   \label{def:model}
    Let $X$ be a smooth and proper scheme over $K$.
     \begin{enumerate} 
      \item  A {\em model of $X$ over $\OO_K$} is an algebraic space $\mathcal{X}$ that is flat and 
         proper over $\Spec\OO_K$, together with an isomorphism $X\isomto \mathcal{X}_K$ from $X$ to the generic fiber of $\mathcal{X}$.
      \item We say that $X$ has {\em good reduction} if there exists a model of $X$ that
         is smooth over $\OO_K$.
    \end{enumerate} 
\end{Definition}

We stress that we have to allow algebraic space models when
dealing with K3 surfaces, which is different from the situation for
curves and Abelian varieties.
In fact, we refer the interested reader to \cite[Section 5.2]{Mat15} for
examples of K3 surfaces with good reduction, where no smooth model
exists in the category of schemes.

\subsection{Inertia and monodromy}
\label{subsec: inertia}
The $G_K$-action on $\overline{K}$ induces an action on $\OO_{\overline{K}}$
and by reduction, an action on $\overline{k}$.
This gives rise to a continuous and surjective homomorphism $G_K\to G_k$ 
of profinite groups. 
Thus, we obtain a short exact sequence
$$
1\,\to\,I_K\,\to\,G_K\,\to\,G_k\,\to\,1,
$$
whose kernel $I_K$ is called the {\em inertia group}.
In fact, $I_K$ is the absolute Galois group 
of the maximal unramified extension $K^{{\rm nr}}$ 
of $K$.
Then, an {\em $\ell$-adic representation of $G_K$} consists
of a finite-dimensional $\QQ_\ell$-vector space $V$ together 
with a continuous group homomorphism $\rho:G_K\to{\rm GL}(V)$.
The representation $\rho$ is called {\em unramified}
if $\rho(I_K)=\{{\rm id}_V\}$.

A relation between good reduction and unramified
Galois representations on $\ell$-adic cohomology groups is given
by the following well-known result, which follows from the proper smooth 
base change theorem.
For schemes, it is stated in \cite[Th\'eor\`eme XII.5.1]{SGA4iii}, 
and in case the model is an algebraic space,
we refer to \cite[Theorem 0.1.1]{LZ14} or \cite[Chapitre VII]{Art73}.

\begin{Theorem}
 \label{thm:obvious}
  Let $X$ be a smooth and proper variety over $K$ with good reduction. 
  Then, the $G_K$-representation
  on $\Het{n}(X_{\overline{K}},\QQ_\ell)$ is unramified for all $n$ and for all 
  $\ell\neq p$.
\end{Theorem}

\subsection{\texorpdfstring{$\bm{p}$}{p}-adic Galois representations}
Now let us assume that $p>0$. In \cite{Fon82,Fon94b}, Fontaine introduced his famous period rings
$$
 \Bcris \,\subseteq\, \Bst \,\subseteq \,\BdR,
$$
which are $\QQ_p$-algebras equipped with various extra structures, the most relevant for us being compatible actions of $G_K$. One usually assumes that $K$ complete, but since completion does not change the absolute Galois group, there is no need to do so. Given a $p$-adic representation $\rho:G_K\to {\rm GL}(V)$,
that is, $V$ is a finite dimensional $\QQ_p$-vector space and $\rho$ is
a continuous group homomorphism, we define
\begin{align*}
  \DDcris (V) &:= \left( V\otimes_{\QQ_p} \Bcris \right)^{G_K} \\
  \DDst (V)   &:= \left( V\otimes_{\QQ_p} \Bst \right)^{G_K} \\
  \DDdR (V) &:= \left( V\otimes_{\QQ_p} \BdR \right)^{G_K} \\
  \DDpst (V) &:= \varinjlim_{K'/K} \left( V\otimes_{\QQ_p} \Bst \right)^{G_{K'}},
\end{align*}
where the limit in the last definition is taken over all finite extensions of $K$ 
(within our fixed algebraic closure $\overline{K}$), and
where all $G_K$-invariants are taken with respect to the diagonal action given by $\rho$
on the first factor and the natural $G_K$-action on the period rings.
Then:
 \begin{enumerate} 
\item $\DDcris(V)$ and $\DDst(V)$ are finite-dimensional $K_0$-vector spaces that come equipped 
 with a semi-linear Frobenius, that is, they are $F$-isocrystals over $k$,
\item $\DDdR(V)$ is a finite dimensional $\widehat{K}$-vector space,
\item $\DDpst(V)$ is a finite-dimensional $K_0^{{\rm nr}}$-vector
space equipped with a semi-linear Frobenius and a semi-linear action of $G_K$.
 \end{enumerate} 
Here, $G_K$-acts on $K_0^{{\rm nr}}$ via its quotient 
$G_k\cong \mathrm{Gal}(K_0^{{\rm nr}}/K_0)$ and in particular, the 
induced $I_K$-action is linear. 
(Of course, this does not exhaust the extra structures these vector spaces have, but we 
will not need to explicitly use either the Hodge filtration or the monodromy operator.) 
By \cite[Propositions 1.4.2 and 5.1.2]{Fon94a}, 
we have inequalities
\begin{align*}
  \dim_{K_0} \DDcris(V) \,\leq\, \dim_{K_0} \DDst(V) \,\leq\,  &\dim_{K_0^{{\rm nr}}} \DDpst(V) \\
 &\,\leq\,\dim_{\widehat{K}}\DDdR(V)\,\leq\, \dim_{\QQ_p} V.
\end{align*}
The representation $\rho$ is called {\em crystalline} (resp. {\em semi-stable}, resp.
{\em potentially semi-stable}, resp. {\em de\thinspace Rham}) 
if we have equality everywhere 
(resp. the last three, two, one inequalities are equalities). We also say that a representation is \emph{potentially crystalline} if it becomes crystalline after a finite extension of $K$.

To obtain results for \emph{integral} $p$-adic cohomology, we will need to recall 
the notion of a Breuil--Kisin module, originally introduced by Kisin in \cite{Kis06}. 
For $\mathfrak{S}:=W\pow{u}$, we have two natural and surjective maps $\mathfrak{S}\to W$ 
and $\mathfrak{S}\to \OO_{\widehat K}$ defined by $u\mapsto 0$ and $u\mapsto \pi$, respectively. 
Let $E(u)$ be the monic Eisenstein polynomial that generates the kernel 
of $\mathfrak{S}\rightarrow \OO_{\widehat K}$. 
There is a Frobenius map $\sigma:\mathfrak{S}\rightarrow \mathfrak{S}$ that is the 
absolute Frobenius on $W$ and sends $u\mapsto u^p$.

\begin{Definition} 
A \emph{Breuil--Kisin module} over $\OO_K$ is a finite
free $\mathfrak{S}$-module $M$ together with a 
morphism of $\mathfrak{S}$-modules
$$
  \varphi\,:\, M \otimes_{(\mathfrak{S},\sigma)}\mathfrak{S} \,\to\, M
$$
that becomes an isomorphism after inverting $E(u)$.
\end{Definition}

Again, by definition this only depends on $\widehat{K}$. By specializing along $u=0$, there is a functor from Breuil--Kisin modules 
over $\cal{O}_K$ to $F$-crystals over $k$. Alternatively, we can consider the closely related functor which specializes along the composite map
\[ \mathfrak{S}\overset{u=0}{\rightarrow} W \overset{\sigma}{\rightarrow} W, \]
where $\sigma$ is the Frobenius map on $W$. 
Then, Kisin showed how to construct a functor $\mathrm{BK}_{\OO_K}$ from 
$G_K$-stable lattices in crystalline $G_K$-representations to Breuil--Kisin modules,
for a precise statement we refer to the reader to \cite[Theorem 1.2.1]{Kis10}. 
In particular, by specialization, we obtain $F$-crystals
$\mathrm{BK}_{\OO_K}(\Lambda)\otimes_{\mathfrak{S}} W$ and $\mathrm{BK}_{\OO_K}(\Lambda)\otimes_{\mathfrak{S},\sigma} W$
associated to such a lattice $\Lambda\subset V$. 
It follows from \cite[Theorem 1.2.1 (1)]{Kis10} that 
there is an isomorphism
$\mathrm{BK}_{\OO_K}(\Lambda)\otimes_{\mathfrak{S}} K_0\cong \DDcris(V)$ 
of $F$-isocrystals.

\subsection{The connection to geometry}
Let us continue to suppose that $p>0$. If $X$ is a smooth and proper variety over $K$, then it is well-known 
(see, for example, \cite[Section 6.3.3]{Ber97a}),
that the $G_K$-representations $\Het{n}(X_{\overline{K}},\QQ_p)$ 
are potentially semi-stable for all $n$.
The $p$-adic analog of Theorem \ref{thm:obvious}
is the following theorem, essentially due to Colmez--Nizio\l{} \cite{CN17} 
and Bhatt--Morrow--Scholze \cite{BMS16}. 
For models of K3 surfaces, which is the case we are interested in, 
a part of it has already been proven by Matsumoto \cite[Section 2]{Mat15}.

\begin{Theorem}
 \label{thm:crysrep}
  Let $X$ be a smooth and proper variety over $K$ and assume that 
  there exists a smooth and proper algebraic space
  $$
     {\cal X}\,\to\,\Spec\OO_K,
  $$ 
  whose generic fiber is $X$ and whose special fiber ${\cal X}_k$ is a scheme.  
  Then,
   \begin{enumerate} 
   \item the $G_K$-representation on $\Het{n}(X_{\overline{K}},\QQ_p)$ is crystalline for all $n$.
   \item For all $n$, there exist isomorphisms
    $$
      \DDcris\left( \Het{n}(X_{\overline{K}},\QQ_p) \right) \,\iso\, \Hcris{n}( {\cal X}_k/K_0)
    $$
    of $F$-isocrystals over $k$.
   \item Assume that $\Hcris{\ast}({\cal X}_k/W)$ is torsion-free.
     Then, also $\Het{\ast}(X_{\overline{K}},\ZZ_p)$ is torsion-free, 
     and for all $n$ there exist isomorphisms
     $$
       \BK_{\OO_K}\left( \Het{n}(X_{\overline{K}},\ZZ_p) \right) \otimes_{\mathfrak{S},\sigma} W 
          \,\iso\, \Hcris{n}\left( {\cal X}_k/W \right)
     $$
     of $F$-crystals over $k$.
   \end{enumerate} 
\end{Theorem}

\prf
We may assume that $K=\widehat{K}$. We first claim that the formal completion $\mathfrak{X}$ of $\cal{X}$ along the special fiber is in fact a formal scheme. 
Indeed, for every $m\geq1$, the reduction of the algebraic space 
$\cal{X}_m:=\cal{X}\times\Spec \OO_K/(\pi^{m+1})$ is equal to 
${\cal X}_k$, which is a scheme.
It follows from \cite[Chapter 3, Corollary 3.6]{Knu71}, that $\cal{X}_m$ is also a scheme, 
hence $\mathfrak{X}=\varinjlim_m \cal{X}_m$ is a formal scheme.

Let $\mathfrak{X}_K$ denote its generic fiber as a rigid space over $K$.
Next, we claim that there is an isomorphism $\mathfrak{X}_K\cong X^\mathrm{an}$. 
(This is well-known for schemes, and is probably also well-known to the experts for algebraic spaces.) 
Choose a presentation
$$
   \cal{X}^1 \rightrightarrows \cal{X}^0 \rightarrow \cal{X} 
$$
of $\cal{X}$ as the coequalizer of an \'etale equivalence relation. 
Let $\mathfrak{X}^i$ denote the formal completion of $\cal{X}^i$ and let $X^i$ be the generic fiber. 
By \cite[Proposition 0.3.5]{Ber96b} we obtain a diagram
$$ 
  \xymatrix{ 
  \mathfrak{X}^1_{K} \ar@<0.5ex>[r] \ar@<-0.5ex>[r] \ar[d] & \mathfrak{X}^0_{K} \ar[r]\ar[d] & \mathfrak{X}_K \\ 
  X^{1,\mathrm{an}} \ar@<0.5ex>[r] \ar@<-0.5ex>[r] & X^{0,\mathrm{an}} \ar[r] & X^\mathrm{an} 
  }
$$
with vertical maps open immersions.
In particular, we get a map $\mathfrak{X}_K\rightarrow X^\mathrm{an}$ by the universal property 
of coequalizers. 
This moreover fits into a commutative diagram
$$ 
 \xymatrix{  
  \mathfrak{X}^0_{K} \ar[r]\ar[d] & \mathfrak{X}_K \ar[d] \\ 
  X^{0,\mathrm{an} }\ar[r] & X^\mathrm{an} 
 }
$$
of rigid analytic spaces over $K$, in which the horizontal arrows are \'etale (of the same degree) 
and the left hand vertical map is an open immersion. 
In particular, $\mathfrak{X}_K\rightarrow X^\mathrm{an}$ is \'etale and since both 
$\mathfrak{X}_K$ and $X^\mathrm{an}$ are proper, it must be finite.
Using again the commutativity of the diagram, we see that it must be of degree one, 
and hence, an isomorphism.

Thus applying \cite[Theorem 3.8.1]{Hub96} we obtain a $G_K$-equivariant isomorphism
$$
   \Het{n}(X_{\overline{K}},\ZZ_p)  \,\cong\, \Het{n}( \mathfrak{X}_{\widehat{\overline{K}}},\ZZ_p ), 
$$
for all $n\geq0$, and the same is also true for $\QQ_p$-valued \'etale cohomology. 
Claims (1) and (2) now follow from 
\cite[Corollary 1.10]{CN17}, and in fact all three follow from \cite[Theorem 1.1]{BMS16}.
\qed\medskip

\begin{Remark} \label{rem: crysrep} 
\begin{enumerate}
 \item \label{num: rem crysrep1} Note that the integral comparison result in \cite{BMS16} is stated without a Frobenius pull-back, however, it was pointed out to us by T. Koshikawa and P. Scholze that this is not quite correct. The problem is that the map $\mathfrak{S}\rightarrow A_\mathrm{inf}$ defined on p.33 of \cite{BMS16} does not make the resulting diagram commute: one needs to compose with the Frobenius map on $W(k)$. Since the comparison isomorphism provided by \cite[Theorem 1.1]{BMS16} is defined by first base changing to $A_\mathrm{inf}$, to get the correct statement over $\mathfrak{S}$ one needs to include this Frobenius pull-back.
 \item As already mentioned in Section \ref{subsec: models},
 one has to work with models that are algebraic spaces.
 The assumption in Theorem  \ref{thm:crysrep} on the special fiber 
 of a smooth and proper model being a scheme is unnatural and it is probably not needed,
 but we had to impose it in order to use the results that are currently available.
 We note that smooth and proper algebraic spaces of dimension $\leq2$ over $k$, 
 as well as smooth and proper algebraic spaces over $k$
 that are group spaces, are schemes by 
 \cite[II.6.7 and V.4.10]{Knu71}. 
 Thus, if ${\cal X}\to\Spec\OO_K$ is a smooth and proper algebraic space 
 model of some variety $X$, such that $X$ is of dimension $\leq2$ or such that
 $X$ is an Abelian variety over $K$, then the 
 special fiber ${\cal X}_k$ is a scheme.
  \end{enumerate}
\end{Remark}

\subsection{A lemma on crystalline representations}
In the proof of Theorem \ref{thm: first main result in introduction}
below, we will need some elementary facts about potentially crystalline and 
potentially semi-stable $p$-adic Galois representations,
which are certainly well-known to the experts. We continue to assume that $p>0$.

\begin{Lemma}
 \label{lemma: change of group}
   Let $\rho:G_K\to {\rm GL}(V)$ be a $p$-adic Galois representation.
   Let $K\subseteq K'$ be a finite field extension and denote by 
   $\rho' :G_{K'}\to {\rm GL}(V)$ the restriction of $\rho$ to $G_{K'}$.
  \begin{enumerate} 
  \item \label{num: lcog1} Assume that $K\subseteq K'$ is unramified.   
   Then, $\rho$ is a crystalline $G_K$-representation if and only if $\rho'$ is a crystalline
   $G_{K'}$-representation.
  \item  Assume that $K\subseteq K'$ is totally ramified. 
   Then, the following are equivalent:
    \begin{enumerate} 
     \item $\rho$ is crystalline,
     \item $\rho'$ is crystalline and the inclusion
        $$
         \left( V\otimes_{\QQ_p} \Bcris \right)^{G_{K'}} \subset  \left( V\otimes_{\QQ_p} \Bcris \right)^{G_{K}}
        $$
        is an equality. 
     \end{enumerate} 
   \item \label{num: lcog3} The following are equivalent:
       \begin{enumerate} 
        \item $\rho$ is crystalline,
        \item $\rho$ is potentially crystalline and the $I_K$-action on $\DDpst(V)$ is trivial.
       \end{enumerate} 
  \end{enumerate} 
\end{Lemma}

\prf
We may assume that $K=\widehat{K}$. First, suppose that $K'/K$ is unramified.
If $\rho$ is crystalline, then so is $\rho'$ and thus, it remains
to show the converse direction.
After replacing the extension $K\subseteq K'$ 
by its Galois closure, we may assume 
that $K'$ is Galois over $K$, say with group $H$.
Then, we have the following invariants with respect to the induced 
$H$-actions
\begin{equation}
 \label{eq: H invariants}
   \left( \left( V\otimes_{\QQ_p} \Bcris \right)^{G_{K'}} \right)^H
   \,\iso\, \left( V\otimes_{\QQ_p} \Bcris \right)^{G_{K}}\,. \tag{$\ast$}
\end{equation}  
Since $\rho'$ is crystalline, 
$\left( V\otimes_{\QQ_p} \Bcris \right)^{G_{K'}}$ is a $K_0'$-vector
space of dimension $d:=\dim_{\QQ_p}(V)$.
Since $K'$ is unramified over $K$, also $K_0\subseteq K_0'$
is Galois with group $H$. 
Thus, the $H$-invariants in \eqref{eq: H invariants} are $K_0$-vector
spaces of dimension $d$.
This proves that $\rho$ is crystalline and establishes 
claim (1).

Next, assume that $K'/K$ is totally ramified.
Then, we have $K'_0=K_0$ and since $G_{K'}$ is contained in $G_K$, 
we find
$$
  \dim_{K_0} \left( V\otimes_{\QQ_p}\Bcris \right)^{G_K}\,\leq\,
  \dim_{K_0} \left( V\otimes_{\QQ_p}\Bcris \right)^{G_{K'}} \,\leq\, \dim_{\QQ_p} V\,.
$$
We have equality on the right if and only if $\rho'$ is crystalline, 
equality on the left if and only if $H$ acts trivially on 
$\left( V\otimes_{\QQ_p}\Bcris \right)^{G_{K'}}$, 
and equality everywhere if and only if $\rho$ is crystalline. 
This establishes claim (2).

For claim (3), we note that since crystalline implies potentially crystalline,
we may assume that $\rho$ is potentially crystalline. 
Next, we choose a finite Galois extension $K'/K$ such that $\rho':=\rho|_{G_{K'}}$ 
is crystalline. 
Applying part (1) and replacing $K$ by its maximal unramified extension inside $K'$,
we may assume that $K'/K$ is totally ramified. 
Denote by $H$ its Galois group. 
Given part (2), it now suffices to observe that in this situation we have
$$
  \DDpst(V) \,\cong\, \left( V\otimes_{\QQ_p}\Bcris \right)^{G_{K'}} \otimes_{K_0} K_0^{{\rm nr}}, 
$$
and the $I_K$-action on $\DDpst(V)$ is simply the extension of scalars of the natural action of 
$I_K$ on 
$$
  \left( V\otimes_{\QQ_p}\Bcris \right)^{G_{K'}}
$$
via the surjective map $I_K\onto H$.
\qed\medskip

\subsection{Forms of $\bm{F}$-isocrystals and non-abelian cohomology} 
\label{subsec: descent}
In this section we recall some material on $L$-forms and descent for $F$-isocrystals. Here, we will assume that $p>0$, and that $K=K_0$ is complete and absolutely unramified.

Let $K'/K$ be a finite and unramified Galois extension with Galois group $G$, and let $V$ be an $F$-isocrystal over $k$. Let $k'$ denote the residue field of $K'$.

\begin{Definition} 
\begin{enumerate} 
  \item A \emph{$k'$-form of $V$} is an $F$-isocrystal $U/K$ such that there exists
     an isomorphism $V\otimes K' \cong U\otimes K'$ as $F$-isocrystals over $k'$. 
  \item Two $k'$-forms of $V$ are said to be \emph{equivalent} if they are isomorphic as 
 $F$-isocrystals over $k$.
  \end{enumerate} 
 We denote the set of equivalence classes of $k'$-forms by $E(k'/k,V)$.
\end{Definition}

We note that $E(k'/k,V)$ is a pointed set, the distinguished element being the class of $V$ itself. 
Given a form $U \in E(k'/k,V)$ then by definition there is an isomorphism
$$
 f \,:\, V\otimes_K K' \,\isomto\, U\otimes_K K'  
$$
of $F$-isocrystals over $k'$. 
Moreover, both sides carry a natural $G$-action coming from the $G$-action on $K'$, 
but in general, $f$ will \emph{not} be equivariant for these actions. 
In fact, we can define a function
\begin{align*} 
  \alpha_f \, :\, G &\rightarrow \mathrm{Aut}_{K',F}\left( V\otimes_K K'\right) \\
  \alpha_f(\sigma) &= f^{-1}\circ \sigma(f)
\end{align*}
that measures the failure of $f$ to be $G$-equivariant. 
We note that $\alpha_f$ is a $1$-cocyle for the natural action of $G$ on 
$\mathrm{Aut}_{K',F}(V\otimes_K K')$. Then, we have the following standard result.

\begin{Proposition} 
 \label{prop: descent} 
 This map induces a bijection
$$ 
 E(k'/k,V)\,\to\, H^1\left(G,\mathrm{Aut}_{K',F}(V\otimes_K K')\right)
 $$
of pointed sets. 
\end{Proposition}

It is worth recalling how to construct the inverse functor, 
via a slightly different perspective on descent. 
That is, as well as the category $F\text{-}\mathrm{Isoc}(k)$ of $F$-isocrystals over $k$, 
we can consider the category $G\text{-}F\text{-}\mathrm{Isoc}(k')$ of $F$-isocrystals over $k'$ together 
with a compatible, semi-linear $G$-action. 
Again, the proof of the following result is straightforward.

\begin{Proposition} 
 \label{prop: better descent} 
 The functors
 $$\begin{array}{lclclcl}
   F\text{-}\mathrm{Isoc}(k) & \rightarrow & G\text{-}F\text{-}\mathrm{Isoc}(k') &:&
   V & \mapsto& V\otimes_K K' \\
   G\text{-}F\text{-}\mathrm{Isoc}(k') & \rightarrow & F\text{-}\mathrm{Isoc}(k) &:&
   W & \mapsto& W^G
  \end{array}$$
 are mutual inverse equivalences of categories.
\end{Proposition}

Using this proposition, the inverse map in Proposition \ref{prop: descent} is then described as follows:
given a $1$-cocycle $\alpha:G\rightarrow \mathrm{Aut}_{K',F}(V\otimes_K K')$, 
we can define a new \emph{action} of $G$ on $V\otimes_K K'$ 
via the \emph{homomorphism}
\begin{align*} 
  \rho^\alpha:G &\rightarrow \mathrm{Aut}_{K',F}(V\otimes_K K')  \\ 
  \rho^\alpha(\sigma)(v) &= \alpha(\sigma)(\sigma(v))
\end{align*}
where $\sigma(v)$ is the ``standard'' action. 
By Proposition \ref{prop: better descent} we obtain an $F$-isocrystal
$$
  V^\alpha \,:=\, (V\otimes_K K')^{\rho^\alpha}
$$
over $k$, which is the form of $V$ giving rise to $\alpha$. 
Note in particular, that $V$ and $V^\alpha$ are isomorphic over $k$ 
if and only if there exists an automorphism
$$
  V\otimes_K K' \,\to\, V\otimes_K K'
$$
of $F$-isocrystals over $k'$ which intertwines the `natural' $G$-action on the left hand side 
with the `$\alpha$-twisted' action on the right hand side.

There is also similar discussion when $K'$ is replaced by 
the maximal unramified extension $K^\mathrm{nr}$, 
as well as for $F$-crystals over $k$. We leave the details to the reader.

\subsection{Forms of $G_k$-modules and non-abelian cohomology} 
\label{subsec: descent2}
There is a similar theory for $\ell$-adic Galois representations
of the absolute Galois group $G_k$, allowing $p\geq 0$ and $\ell$ to be any prime, including possibly $p$.

If we let $k'/k$ be a finite Galois extension, with Galois group $G$, then we have an exact sequence
$$
1 \,\to\, G_{k'} \,\to\, G_k \,\to\, G \,\to\, 1, 
$$
where $G_{k'}$ is the absolute Galois group of $k'$. 

\begin{Definition} 
 Let $V$ be a finite dimensional $\QQ_\ell$ vector space 
 and $\rho:G_k\to\mathrm{GL}(V)$ be an $\ell$-adic Galois representation.
  \begin{enumerate} 
  \item A \emph{$k'$-form of $V$} is an $\ell$-adic Galois representation
    $\psi:G_k\to\mathrm{GL}(W) $ such that there exists
    a $G_{k'}$-equivariant and $\QQ_\ell$-linear isomorphism 
    $V|_{G_{k'}} \cong W|_{G_{k'}}$.
  \item Two $k'$-forms of $V$ are said to be \emph{equivalent} if they are isomorphic as 
   $G_k$-representations.
  \end{enumerate} 
We denote the set of equivalence classes of $k'$-forms by $E(k'/k,V)$. 
\end{Definition}

As in the previous section, $E(k'/k,V)$ is a pointed set. 
There is a natural $G_k$ action on $\mathrm{Aut}_{G_{k'}}(V)$ via conjugation, and by definition $G_{k'}$ acts trivially. 
Thus, this action factors through $G$. Given a $k'$-form $W\in E(k'/k,V)$, we can choose a $G_{k'}$-equivariant isomorphism
$$ 
f \,:\, V \,\to\,  W
$$
and we obtain a $1$-cocycle
\begin{align*} 
  \alpha_f \, :\, G &\rightarrow \mathrm{Aut}_{G_{k'}}(V) \\
  \alpha_f(\sigma) &= f^{-1}\circ \sigma(f)
\end{align*}
as before. 
Then, we have the following proposition, whose proof we leave to the reader.

\begin{Proposition} 
 \label{prop: descent2} 
This map induces a bijection
 $$ 
   E(k'/k,V)\,\to\, H^1\left(G,\mathrm{Aut}_{G_{k'}}(V)\right)
 $$
of pointed sets.
\end{Proposition}

We can describe the inverse functor relatively easily as follows. 
Namely, if we have a $G_k$-representation $\rho:G_k \rightarrow \mathrm{GL}(V)$ and a $1$-cocycle
$$
\alpha\,:\,G \,\to\, \mathrm{Aut}_{G_{k'}}(V),
$$
then we can construct a ``new'' $G_{k}$-action on $V$ by defining
$$
\rho^\alpha(g)(v) \,:=\, \alpha(g)(\rho(g)(v))
$$
Then, $\rho^\alpha:G_k\rightarrow \mathrm{GL}(V)$ is the $k'$-form of $V$ 
giving rise to $\alpha$. All of this works integrally as well, that is, replacing $\QQ_\ell$ by $\ZZ_\ell$, and there is a version replacing $k'/k$ with $\bar k/k$. We leave all these details up to the reader.

\section{Canonical surface singularities and Weyl groups} \label{sec: root}
\label{subsec: RDP and K3}
We will need to fix a couple of definitions concerning surfaces and their
singularities, as well as make a few group-theoretic calculations involving Dynkin diagrams and Weyl groups.

\begin{Definition}
  \label{def: singular surfaces}
  Let $F$ be a field.
   \begin{enumerate} 
   \item A {\em surface} over $F$ is a separated, geometrically integral
     scheme $S$ of finite type and dimension $2$ over $\Spec F$.
    \item A surface $S$ over $F$ has {\em (at worst) canonical singularities}
     if it is geometrically normal, Gorenstein, and if the 
     minimal resolution of singularities $f:\widetilde{S}\to S$
     satisfies $f^\ast\omega_S\cong\omega_{\widetilde{S}}$.
     Here, $\omega_S$ and $\omega_{\widetilde{S}}$ denote the
     dualizing sheaves of $S$ and $\widetilde{S}$, respectively.
   \item A {\em K3 surface} over $F$ is smooth and proper surface $S$ over 
     $F$ with $\omega_S\cong\OO_S$ and $H^1(S,\OO_S)=0$.
   \end{enumerate} 
\end{Definition}

We note that canonical surface singularities are rational singularities,
and that canonical surface singularities are characterized as being those
rational singularities that are Gorenstein.
If the ground field $F$ is algebraically closed, then canonical surface 
singularities are also known as {\em rational double point singularities}, 
{\em du Val singularities}, or  {\em Kleinian singularities}, 
and we refer to Artin's articles \cite{Art62,Art66} 
for details and proofs.

We also note here that an algebraic space that is
two-dimensional, smooth, and proper over a field is a projective scheme
\cite[V.4.10]{Knu71}.
By \cite[Theorem 2.3]{Art62}, the same is true for two-dimensional
and proper algebraic spaces with rational singularities, which
includes canonical surface singularities.

\subsection{Classification of canonical singularities}

If $F$ is a field, then canonical surface singularities over $F$ are classified 
by the finite Dynkin diagrams, as explained in \cite[\S24]{Lip69} or \cite[\S3.3]{Kol13}. 
Let us recall here a little about how this works, as well as fixing some notations that will be used later on.

Let $S/F$ be a surface with canonical singularities and let $\widetilde{S}\rightarrow S$ be its minimal resolution. 
Then, the possible configurations of the exceptional divisor $E_P\subset \widetilde{S}$ over some singular point 
$P\in S$ are described by the Dynkin diagrams
$$
 A_n, \;B_n, \;C_n, \;D_n, \;E_6, \;E_7, \;E_8, \;F_4, \;G_2
 $$
(whenever we say ``Dynkin diagram'', we will always mean a \emph{finite, reduced, irreducible} Dynkin diagram, 
that is, one from the above list). 
The geometric interpretation of a diagram as follows. 
Nodes of the diagram correspond to irreducible (but not necessarily geometrically irreducible) components $E_{P,i}\subset E_{P}$, and we have
$$
 n_{P,i} = h^0(E_{P,i}) \,:=\, \dim_{F(P)}H^0(E_{P,i},\cal{O}_{E_{P,i}}),
$$
where $n_{P,i}$ is the integer labelling the node corresponding to $E_{P_i}$, and where $F(P)$ is the residue field at $P$. 
The existence of an edge between $E_{P,i}$ and $E_{P,j}$ means that $E_{P,i}\cap E_{P,j}\neq \emptyset$. 
This then completely determines the intersection matrix of the $E_{P_i}$ by 
$$
E_{P,i}\cdot E_{P,j} \,=\, 
\begin{cases} 
  -2 n_{P,i}  &\text{ if }i=j \\ 
  \mathrm{max}(n_{P,i},n_{P,j})=\mathrm{lcm}(n_{P,i},n_{P,j})&\text{ if } i\neq j\text{ and }E_{P,i}\cap E_{P,j}\neq \emptyset \\ 
  0 &\text{ if } E_{P,i}\cap E_{P,j}= \emptyset 
\end{cases} 
$$
(see for example \cite[Lemma 22.2]{Lip69} and its proof). Of course, we have $E_{P,i}\cdot E_{Q,j}=0$ if $P\neq Q$. 

We will call the subgroup $\Lambda_S \subseteq \mathrm{Pic}(\widetilde{S})$ that is (freely) generated by the integral components 
of all $E_P$ (as $P$ ranges over the singular points of $S$), together with its natural intersection pairing, 
the \emph{root lattice} of $S$, and its Weyl group
$$
\cal{W}_S \,\leq\, \mathrm{Aut}_{\ZZ}(\Lambda_S) 
$$
the \emph{Weyl group} of $S$. 
Concretely, $\cal{W}_S$ is generated by the maps
$$ 
\begin{array}{ccccc}
  s_{P,i} &:& \Lambda_S & \to&\Lambda_S\\
    & & D &\mapsto& D\,+\,\frac{1}{n_{P,i}}(D\cdot E_{P,i})\,\cdot\,E_{P,i}
 \end{array}
$$
for all integral components $E_{P,i}$ of all $E_P$. 
Using the same formula, we can view $\cal{W}_S$ as a subgroup of $\mathrm{Aut}_{\ZZ}(\mathrm{Pic}(\widetilde{S}))$.

Note that the natural ``geometric'' sign convention we are using is the opposite of that usually found 
when discussing root systems, for us therefore Cartan matrices will be negative definite rather than positive definite. 
Also note that the root lattice and the Weyl group break up into direct products of the root lattices and Weyl groups associated 
to each individual singular point of $S$.

It is important to remember that these definitions are not stable under a base field extension $F'/F$, 
\emph{even if all the singularities are $F$-rational and of type ADE}. 
For example, an $F$-rational singularity of type $A_1$ could become a singularity of type $A_2$ geometrically.

Nonetheless, we can describe their behaviour under such extensions as follows, at least when 
the extension $F'/F$ is Galois, say with group $G$. 
Indeed, in this case there are natural and compatible $G$-actions on both 
$\Lambda_{S_{F'}}$ and $\cal{W}_{S_{F'}}$, and we clearly have $\Lambda_S = \Lambda_{S_{F'}}^G$. 
Moreover, we can construct a homomorphism 
$$
\cal{W}_S \,\to\, \cal{W}_{S_{F'}}^G 
$$
fom $\cal{W}_S$ into the $G$-invariants of $\cal{W}_{S_{F'}}$ as follows. Take a reflection $s_{P,i}$ in the class of an integral component $E_{P,i}$ of the exceptional divisor $E$, and suppose that $E_{P,i}$ splits into components $\{E_{Q,1},\ldots,E_{Q,n}\}$ over $F'$. 
Then, 
 \begin{enumerate}  
 \item either $E_{Q,j}\cdot E_{Q,j'}=0$ for all $j\neq j'$,
 \item or else $n$ is even and, after possibly re-ordering, we have 
   $E_{Q,2j-1}\cdot E_{Q,2j}=1$ for all $j$ and $E_{Q,j}\cdot E_{Q,j'}=0$ for all other $j\neq j'$,
  \end{enumerate} 
where we argue as in the proof of \cite[Lemma 4.3]{LM18} to see that these are the only cases. In the first case, we simply map $s_{P,i}$ to the product $\prod_j s_{Q,j}$ of the reflection in each of the 
$E_{Q,j}$ (in any order, since they all commute). 
In the second case, we map $s_{P,i}$ to
$$
   \prod_{j=1}^{\frac{n}{2}} s_{Q,2j-1}\cdot s_{Q,2j}\cdot s_{Q,2j-1} = \prod_{j=1}^{\frac{n}{2}}s_{Q,2j}\cdot s_{Q,2j-1}\cdot s_{Q,2j}.
$$
We wish to prove that this induces an isomorphism 
$$
 \cal{W}_{S}\,\isomto\,\cal{W}_{S_{F'}}^G ,
$$
and this amounts to a lemma on invariants of Weyl groups under folding, that surely must be well-known to the experts. 

So suppose that we have a Dynkin diagram $T$, and a group $G \leq \mathrm{Aut}(T)$ of automorphisms of $T$. 
As before, we can see the root lattice $\Lambda_{T/G}$ of $T/G$ as the $G$-invariants $\Lambda_T^G$ 
of the root lattice of $T$, and exactly as above we may construct a homomorphism
$$
   \cal{W}_{T/G} \,\rightarrow\, \cal{W}_T^G 
 $$
from the Weyl group of the quotient diagram $T/G$ to the $G$-invariants of the Weyl group of $T$. 

\begin{Lemma} 
 \label{lemma: GinvW} 
 This homomorphism $ \cal{W}_{T/G} \rightarrow \cal{W}_T^G$ is an isomorphism.
\end{Lemma}

\prf 
We split into two cases, depending on whether each orbit of $G$ on the simple roots in $\Lambda_T$ 
only contains pairwise orthogonal roots or not.

If the simple roots in each orbit are pairwise orthogonal, then we use an argument we learned from \cite{Ste08}. 
It is easy to verify that if $w \in \cal{W}_{T/G}$ maps to $\bar w \in \cal{W}_T$, 
then $w$ and $\bar w$ agree as endomorphisms of 
$\Lambda_{T/G}=\Lambda_T^G\subset \Lambda_T$, thus the given map is injective. 
To see surjectivity, that is, that every $w \in \cal{W}_T^G$ is in fact an element of $\cal{W}_{T/G}$ we induct on 
the \emph{length} of $w$ in the sense of \cite[\S1.6]{Hum90}. 
If $\ell(w)=0$, then there is nothing to prove, and if $\ell(w)>0$ then we can find a basic reflection $s_\alpha$ 
in some simple root $\alpha\in \Lambda_T$ such that $\ell(ws_\alpha)<\ell(w)$. Thus $\ell(ws_\alpha) =\ell(w)-1$ and $w(\alpha)<0$ by \cite[Lemma 1.6, Corollary 1.7]{Hum90}.

Let $O(\alpha)$ be the orbit of $\alpha$ under the $G$-action, and $s:=\prod_{\beta\in O(\alpha)} s_\beta$. Since any automorphism of $T$ is length preserving we see that we also have 
$\ell(ws_{\beta})=\ell(w)-1$ and $w(\beta)<0$ for any $\beta\in O(\alpha)$. Now let $w'=ws$, we claim that $\ell(w')<\ell(w's_\beta)$ for all $\beta \in O(\alpha)$. Indeed, since the $\beta$ are pairwise orthogonal, we find that $w'(\beta)=w(-\beta)>0$ and thus again that $\ell(w')=\ell(w's_\beta)-1$.

Therefore $w'$ is a minimal coset representative for the parabolic subgroup generated by the $\{s_\beta\; |\; \beta\in O(\alpha)\}$ in the sense of \cite[\S1.10]{Hum90}. Since $w=w's$ we can therefore deduce that $\ell(w)=\ell(w')+\ell(s)$, and again since the $\beta$ are pairwise orthogonal that $\ell(s)=\lvert O(\alpha) \rvert$. Thus, $\ell(w')<\ell(w)$, and since $s$ is clearly $G$-invariant, so is $w'$. We may therefore apply the induction hypothesis to conclude.

In the second case, the only possibility is that $T=A_{2n}$, $G=\ZZ/2\ZZ$ and $T/G=C_n$. 
If we label the vertices of $A_{2n}$ in the obvious way, then we have $\cal{W}_T\cong S_{2n+1}$, 
mapping the reflection $s_{\alpha_i}$ in the $i$th basis element to the transposition $(i,i+1)$, see \cite[\S2.10]{Hum90}. Moreover, the action of $\ZZ/2\ZZ$ on $\mathcal{W}_T\cong S_{2n+1}$ is via conjugation by the element
$$
\sigma \,:=\, (1,2n+1)(2,2n)(3,2n-1)\ldots(n,n+2). 
$$
Also by \cite[\S2.10]{Hum90}, we know that $\mathcal{W}_{T/G} \cong (\ZZ/2\ZZ)^n \rtimes S_n$, via the obvious action of $S_n$ on $(\ZZ/2\ZZ)^n$, and the map $\mathcal{W}_{T/G}\rightarrow \mathcal{W}_{T}$ corresponds to the map sending $(i,i+1)$ to $(i,i+1)(2n+1-i,2n+2-i)$, and the $n$th canonical basis element of $(\ZZ/2\ZZ)^n$ to $(n,n+2)$. One easily verifies that this map is injective, with image precisely the centraliser of $\sigma$.
\qed\medskip

\begin{Corollary} 
 \label{cor: GinvW} 
 If $S$ is a surface over $F$ that has canonical singularities
 and $F'/F$ is a finite Galois extension with group $G$,  
 then the above map
$$  
   \cal{W}_{S} \,\rightarrow\, \cal{W}_{S_{F'}}^G
 $$
 is an isomorphism.
\end{Corollary}

\begin{Remark} 
  \label{rem: weylext} 
  It is important to note that the extensions of elements of $\cal{W}_S$ to $\cal{W}_{S_{F'}}$ are \emph{not} generally defined 
  by the naive extension of the formula for the basic reflections. 
  As we saw, if $E_{P,i}$ is an integral component of the exceptional locus of $\widetilde{S}\rightarrow S$ 
  such that $E_{P,i}$ breaks into two connected components $E_{Q_1,i}$ and $E_{Q_2,i}$ over $F'$, 
  then the extension of $s_{P,i}$ to $\mathrm{Aut}_{\ZZ}(\Lambda_{S_{F'}})$ provided by Corollary \ref{cor: GinvW}
  is given by $s_{Q_1,i}s_{Q_2,i}$. 
  This only agrees with the ``naive'' extension
  $$
      D \,\mapsto\, D\,+\,\frac{1}{n_{P,i}}(D\cdot E_{P,i})\,\cdot\, E_{P,i}
  $$
on the subspace of $\Lambda_{S_{F'}}$ fixed by $G$, that is, $\Lambda_{S}$. 
  More generally, these two formulae agree on the subspace 
  $\mathrm{Pic}(\widetilde{S})\subseteq \mathrm{Pic}(\widetilde{S}_{F'})$.
\end{Remark}

We will also need the following rather ugly lemma, that will be used in the proof of Theorem \ref{theo: flopref}.

\begin{Lemma} 
 \label{lemma: dynkaut} 
 Let $A$ be a (negative definite) Cartan matrix of ADE type, and $C$ the permutation matrix coming from an automorphism 
 of the associated Dynkin diagram. If $C\neq I$, then the matrix $(C-I)A^{-1}$ does not have all integral entries.
\end{Lemma}

\prf 
We proceed on a case by case basis, in cases $E_7$ and $E_8$ there is nothing to prove. 
For $A_n$, we order the vertices in the obvious way. 
Then, the only non-trivial possibility for $C$ is the matrix with $1$'s along the anti-diagonal. 
Therefore, the top row of $(C-I)$ is $(-1,0,0,\ldots,0,1)$ and hence the top left hand entry is
$$
\left((C-I)A^{-1}\right)_{1,1} \,=\,(A^{-1})_{n,1}-(A^{-1})_{1,1}.
 $$
The inverses of the Cartan matrices are computed in \cite[p.95]{Ros97}, and hence,
 by our conventions, we have $(A^{-1})_{n,1}=-1/(n+1)$ and $(A^{-1})_{1,1}=-n/(n+1)$. 
So $\left((C-I)A^{-1}\right)_{1,1}$ is not integral (unless $n=1$). 
Similarly, in the $D_n$ case with $C$ of order $2$, we order the vertices as follows.
\begin{figure}[H]\begin{tikzpicture}
\tikzstyle{every node}=[draw,circle,fill=black,minimum size=5pt,inner sep=0pt]
 \node[label={[label distance=1mm]150: \small{$n-2$}}](1) at (0,0) {};
 \node[](2)   at (-1.5,0) {};   
 \node[label={[label distance=2mm]0: \small{$n-1$}}](3) at (0.75,1.3) {};
 \node[label={[label distance=2mm]0: \small{$n$}}](4)   at (0.75,-1.3)  {} ;
 \node[](5) at (-3,0) {};
 \node[label={[label distance=2mm]90: \small{$2$}}](6) at (-4.5,0)  {};
  \node[label={[label distance=2mm]90: \small{$1$}}](7) at (-6,0)   {};

 \path[line width=.3mm] (1) edge (2);
  \path[line width=.3mm] (1) edge (3);
   \path[line width=.3mm] (1) edge (4);
      \path[line width=.3mm] (2) edge (5);
         \path[line width=.3mm] (5) edge (6);
         \path[line width=.3mm] (6) edge (7);
\end{tikzpicture}\end{figure}
\noindent Then, possibly up to reordering when $n=4$, $C$ is the matrix
$$
\left( \begin{matrix} I_{n-2} &  & \\ & 0 & 1 \\ & 1 & 0 \end{matrix}  \right)  
$$
and hence,
$$  
\left((C-I)A^{-1}\right)_{n,n} \,=\, (A^{-1})_{n-1,n}-(A^{-1})_{n,n}.
$$
We again appeal to \cite[p.95]{Ros97} (remembering the change in sign convention) to deduce that 
the bottom right entry $\left((C-I)A^{-1}\right)_{n,n}=1/2$ is not integral.

The remaining cases are handled similarly, by explicitly writing down $C$ and the inverses of the Cartan matrices associated to $D_4$ and $E_6$. In the $D_4$ case with $C$ of order $3$, one finds that 
\[ \left((C-I)A^{-1}\right)_{1,1} \,=\, 1/2  \]
is not integral. In the $E_6$ case with the labelling
\begin{figure}[H]\begin{tikzpicture}
\tikzstyle{every node}=[draw,circle,fill=black,minimum size=5pt,inner sep=0pt]

\node[label={[label distance=2mm]150: \small{$6$}}](3) at (-3,1.5) {};

 \node[label={[label distance=2mm]150: \small{$5$}}](1) at (0,0) {};
 \node[label={[label distance=2mm]150: \small{$4$}}](2)   at (-1.5,0) {};   

 \node[label={[label distance=2mm]150: \small{$3$}}](5) at (-3,0) {};
 \node[label={[label distance=2mm]150: \small{$2$}}](6) at (-4.5,0)  {};
  \node[label={[label distance=2mm]150: \small{$1$}}](7) at (-6,0)   {};

 \path[line width=.3mm] (5) edge (3);
\path[line width=.3mm] (1) edge (2);
      \path[line width=.3mm] (2) edge (5);
         \path[line width=.3mm] (5) edge (6);
         \path[line width=.3mm] (6) edge (7);
\end{tikzpicture}\end{figure}
\noindent and with $C$ necessarily of order $2$, one finds that
\[ \left((C-I)A^{-1}\right)_{1,1} \,=\, 1/3 \]
is not integral. 
\qed\medskip

\begin{Remark} 
 \label{rem: Dn order} 
 We will refer to the ordering of the vertices of $D_n$ described during the proof as the ``usual'' ordering.
\end{Remark}

\section{Some calculations in group cohomology} \label{sec: nagc}

The proof of our N\'eron--Ogg--Shafarevich criterion for K3 surfaces (Theorem \ref{theo: intro NOS}) will essentially depend on showing 
that a certain map in non-abelian cohomology has trivial kernel. 
This can be reduced to a completely group theoretic calculation 
involving the Weyl groups of Dynkin diagrams, and the purpose of
this section is to perform this calculation.

So let $\mathbf{T}=(T_1,\ldots,T_d)$ be a finite collection of 
Dynkin diagrams (that is, each $T_i$ is finite, reduced and irreducible), 
with corresponding root lattice 
$\Lambda_{\mathbf{T}}=\bigoplus_{i=1}^d \Lambda_{T_i}$ 
and Weyl group $\cal{W}_\mathbf{T}=\prod_{i=1}^d\cal{W}_{T_i}$. 
Let $G$ be a finite group acting on $\mathbf{T}$, and let $F$ be a field of characteristic $0$.
\begin{Theorem} 
 \label{theo: nagc main} 
  The induced map
 $$
    H^1\left(G,\cal{W}_\mathbf{T}\right) \,\to\, H^1\left(G,\mathrm{GL}(\Lambda_{\mathbf{T},F})\right)
  $$
 in non-abelian cohomology has trivial kernel.
\end{Theorem}

\begin{Remark} 

\begin{enumerate}
\item Note that this implies that the corresponding result holds with $F$ replaced by any torsion free integral domain, for example, the integers $\ZZ$. 
\item Since the $G$-modules under consideration are non-abelian, 
  $H^1$ is just a pointed set, not a group. 
  But it still makes sense to speak of the kernel of a map of pointed sets $(A,*)\rightarrow (B,*)$, that is, the fiber over the distinguished point of $B$. 
  Since this kernel always contains the distinguished point of $A$, 
  triviality then means that it \emph{only} contains this point. 
  Note that this is a weaker claim than the map being injective, and in fact in the 
  situation of Theorem \ref{theo: nagc main}, the map
$$
     H^1\left(G,\cal{W}_\mathbf{T}\right) \,\to\, H^1\left(G,\mathrm{GL}(\Lambda_{\mathbf{T},F})\right)
$$
will \emph{not} be injective in general, for example it is not injective for the faithful $\ZZ/2\ZZ$ action on $D_5$ (see the proof of Proposition \ref{prop: nagc faith} below).
\end{enumerate}
\end{Remark}

Before we start on the proof of Theorem \ref{theo: nagc main} properly, we will need a lemma telling us how to compute group cohomology in a particular case.

\begin{Lemma} \label{lemma: H1cyc} Let $H$ be a group and $h\in H[n]$ an $n$-torsion element. Let $G=\langle \left. g\right\vert g^n=1\rangle$ be a cyclic group of order $n$, acting on $H$ via conjugation by $h$, and let $H_0 \subset H$ be a subgroup stable under $G$. Then the `twisted' evaluation map 
\begin{align*}
 Z^1(G,H_0) &\rightarrow H \\
 \alpha &\mapsto \alpha(g)h
\end{align*}
induces an injection
\[ H^1(G,H_0) \rightarrow H[n]/\sim \]
from the non-abelian cohomology set $H^1(G,H_0)$ to the set of $H_0$-conjugacy classes of elements of order $n$ in $H$. The image consists of classes meeting the coset $H_0h$.
\end{Lemma}

\begin{proof}
First of all, the cocycle condition tells us that
\begin{align*} 
1=\alpha(1)=\alpha(g^n) &=\alpha(g)\alpha(g)^g\alpha(g)^{g^2}\ldots\alpha(g)^{g^{n-1}} \\ &= \alpha(g)h\alpha(g)h\ldots \alpha(g)h^{-n+1} \\ &= (\alpha(g)h)^n  \end{align*}
and hence indeed $\alpha(g)h \in H[n]$. Next, $\alpha$ and $\tilde{\alpha}$ are cohomologous if, and only if, there exists some $h_0\in H_0$ such that 
\[ h_0\tilde{\alpha}(g) = \alpha(g)hh_0h^{-1} \Leftrightarrow \tilde{\alpha}(g)h = h_0^{-1}\alpha(g)hh_0,  \]
and so the induced map
\[ H^1(G,H_0) \rightarrow H[n]/\sim \]
is injective as claimed. Finally, if $h'\in H[n] \cap H_0h$ then the function
\begin{align*}
\alpha: G\rightarrow H_0\\
\alpha(g^i) = h'^{i}h^{-i}
\end{align*}
is a cocycle such that $\alpha(g)h=h'$ ($\alpha$ takes values in $H_0$ by the hypothesis that $H_0$ is stable under conjugation by $h$).
\end{proof}

The hardest part of the proof of Theorem \ref{theo: nagc main} turns out to be the following special case.

\begin{Proposition} 
 \label{prop: nagc faith} 
 In the situation of Theorem \ref{theo: nagc main} assume that $G$ is cyclic, 
 $\mathbf{T}=T$ is a single Dynkin diagram, and that the action is faithful, that is, $G\rightarrow \mathrm{Aut}(T)$ is injective. 
 Let $\Lambda$ be the root lattice of $T$ and $\cal{W}$ the Weyl group. 
 Then, the induced map
$$
H^1\left(G,\cal{W}\right) \,\to\, H^1\left(G,\mathrm{GL}(\Lambda_{F})\right)
 $$
has trivial kernel.
\end{Proposition}

\prf
We note that if $G$ is non-trivial, then $T$ is simply-laced, that is, it is of type ADE. 

First, assume that $G=\ZZ/2\ZZ$ and $T$ is of type 
$A_n$ with $n\geq 2$, of type $D_n$ with $n\geq 5$ odd, or of type $E_6$. 
Let $g$ denote the unique non-trivial element of $G=\mathrm{Aut}(T)$. Applying \cite[Theorem 1]{KM83}
we know that $-I\in \mathrm{O}(\Lambda) = \cal{W} \rtimes G$, and since $-I\notin \cal{W}$ we can conclude that $g=-w_0$ for some $w_0\in \cal{W}$. 

In particular, the action of $g\in \ZZ/2\ZZ$ on both $\cal{W}$ and $\mathrm{GL}(\Lambda_{F})$ is via conjugation by $w_0$. 
By Lemma \ref{lemma: H1cyc}, we therefore have two bijections
$$
 \begin{array}{lcllcl}
  H^1\left(\ZZ/2\ZZ,\cal{W}\right) &\rightarrow& \cal{W}[2] / \sim, &
  \alpha &\mapsto& \alpha(g)w_0 \\
  H^1\left(\ZZ/2\ZZ,\mathrm{GL}(\Lambda_{F})\right) &\rightarrow& \mathrm{GL}(\Lambda_{F})[2] / \sim, &
  \alpha &\mapsto& \alpha(g)w_0 \\
 \end{array}
$$
between the relevant cohomology sets and the sets of conjugacy classes 
of elements of order $2$. 
We thus need to prove that if $w\in \cal{W}[2]$ is such that $w$ and $w_0$ 
are conjugate inside $\mathrm{GL}(\Lambda_{F})$, then they are already 
conjugate inside $\cal{W}$. 

The conjugacy classes of involutions in Weyl groups are classified 
in \cite[\S4, Theorem A']{Ric82}: they correspond to $\cal{W}$-orbits 
of sub-diagrams $T'\subseteq T$ (where $T'$ is no longer assumed to be irreducible) 
that contain $-I$
in their Weyl group. 
Moreover, the size of $T'$ determines the conjugacy class inside 
$\mathrm{GL}(\Lambda_{F})$, since it determines the multiplicity of $-1$ as an eigenvalue. 
Hence, it suffices to prove that there is a unique orbit of maximal size, 
and that $w_0$ is an element of the corresponding conjugacy class. (Note also that this argument, combined with the examples given in \cite[\S4.1]{Ric82}, shows the failure of injectivity of $H^1\left(G,\cal{W}\right) \to H^1\left(G,\mathrm{GL}(\Lambda_{F})\right)$ when $T=D_5$, since there are at least two distinct $\mathcal{W}$-orbits of such $T'\subset D_5$ of type $2A_1$.)

For $A_n$, a maximal such sub-diagram $T'$ is given by starting at one end 
and choosing alternating vertices. 
There is a unique such $T'$ if $n$ is odd, and two possible choices when $n$ is even. 
Since $w_0=-g$, we see that these two choices are conjugate under $\cal{W}$.  
For $D_n$, with $n\geq 5$ odd, a maximal such diagram is provided 
by $D_{n-1}\subset D_n$, and for $E_6$ by $D_4\subset E_6$. 
In both cases, there is a unique choice for $T'$. 
To conclude, then, we can simply calculate the eigenvalues of $w_0=-g$ 
and observe that we get the maximal possible multiplicity of $-1$ in all 
cases.

The next case we consider is $\ZZ/2\ZZ$ acting on $D_n$ for $n\geq 4$ even 
(although we don't really need to assume $n$ even in the proof). 
In this case by choosing a suitable basis for both $\Lambda_{D_n,F}$ and $\Lambda_{B_n,F}$ we can use the description of Weyl groups given in \cite[\S2.10]{Hum90} to construct a commutative diagram
$$ 
\xymatrix{  (\ZZ/2\ZZ)^{n-1} \rtimes S_n \ar@{^(->}[d]\ar[r]^-{\cong} & \cal{W}_{D_n} \ar[r]\ar@{^(->}[d] & \mathrm{GL}_n(F) \ar@{=}[d] \\
 (\ZZ/2\ZZ)^n \rtimes S_n \ar[r]^-{\cong} & \cal{W}_{B_n} \ar[r] & \mathrm{GL}_n(F)   } 
$$
such that the non-trivial element of $\ZZ/2\ZZ$ acts on both $\cal{W}_{D_n}$ 
and $\mathrm{GL}_n(F)$ as conjugation by the element $w_0\in \cal{W}_{B_n}$ 
corresponding to the element $\delta_n=((0,0,\ldots,0,1),e)\in (\ZZ/2\ZZ)^n \rtimes S_n$. 
Note that the inclusion
$$
   (\ZZ/2\ZZ)^{n-1} \rtimes S_n \hookrightarrow  (\ZZ/2\ZZ)^{n} \rtimes S_n 
$$
on the first factor is as the subspace of elements mapping to zero under the summation map
\[ \sum: (\ZZ/2\ZZ)^{n} \rightarrow \ZZ/2\ZZ.\] 
We first claim that the induced map
$$
  H^1\left(\ZZ/2\ZZ,\cal{W}_{D_n}\right) \,\to\, H^1\left(\ZZ/2\ZZ,\cal{W}_{B_n}\right)
$$
has trivial kernel. As above, using Lemma \ref{lemma: H1cyc} this boils down to showing that if we take an element $x\in  (\ZZ/2\ZZ)^{n-1}\rtimes S_n$ such that $x\delta_n$ and $\delta_n$ 
are conjugate inside $(\ZZ/2\ZZ)^{n} \rtimes S_n$, then they are 
conjugate by an element of $(\ZZ/2\ZZ)^{n-1} \rtimes S_n$. 
This follows from the formula
$$
(\delta,\sigma)^{-1}\,(\delta_n,e)\,(\delta,\sigma) \,=\, (\sigma^{-1}(\delta_n),e) 
$$
for any $(\delta,\sigma)\in (\ZZ/2\ZZ)^{n} \rtimes S_n$, since the RHS does not depend on $\delta$.

Given this claim, it therefore suffices to show that
$$ 
  H^1\left(\ZZ/2\ZZ,\cal{W}_{B_n}\right) \,\to\, H^1\left(\ZZ/2\ZZ,\mathrm{GL}(\Lambda_{B_n,F}) \right)
$$
also has trivial kernel. Since the action is conjugation by $w_0$, as before this amounts to showing that if 
$w\in \cal{W}_{B_n}$ is an element of order $2$ that is conjugate to $w_0$ 
inside $\mathrm{GL}(\Lambda_{B_n,F})$, then it is already conjugate to $w_0$ 
inside $\cal{W}_{B_n}$. 
Again, since the element $w_0$ has $-1$ as an eigenvalue with multiplicity one
and all other eigenvalues $+1$, it suffices to observe by \cite[\S4, Theorem A']{Ric82} 
that all sub-diagrams of $B_n$ of size $1$ are conjugate 
under $\cal{W}_{B_n}$ (see for example the discussion in \S4.1 of \emph{loc. cit.}).

Finally, if we look at $\ZZ/3\ZZ$ acting on $D_4$, then it seems more difficult to give 
such a straightforward description of the action in terms of elements of some Weyl group. Instead we will calculate the kernel of
$$
H^1\left(\ZZ/3\ZZ,\cal{W}\right) \,\to\, H^1\left(\ZZ/3\ZZ,\mathrm{GL}(\Lambda_{F})\right)
$$
by brute force as follows: we order the nodes of $D_4$ in the ``usual'' way (Remark \ref{rem: Dn order}) 
and denote the reflection in the $i$th node by $s_i$. 
Let $\alpha_i$ be the corresponding basis elements for $\Lambda$. 
Then, a generator $g$ for $\ZZ/3\ZZ$ acts via
\begin{align*} 
  \alpha_1 \mapsto \alpha_3,&\;\; \alpha_2 \mapsto \alpha_2 \\
  \alpha_3\mapsto \alpha_4,&\;\;\alpha_4\mapsto \alpha_1.
\end{align*}
By Lemma \ref{lemma: H1cyc}, the kernel of 
$$
H^1\left(G,\cal{W}\right) \,\to\, H^1(G,\mathrm{GL}\left(\Lambda_{F})\right)
$$
is in bijection with $\mathcal{W}$-conjugacy classes of elements of the coset $\mathcal{W}g$ which are conjugate to $g$ in $\mathrm{GL}(\Lambda_F)$. The \cite{SAGE} code reproduced in Appendix \ref{app: code D4} calculates 
an upper bound for the number of elements in this kernel. 
It does so by first listing all elements of $\cal{W}$ in matrix form.
It then calculates which elements $w$ give rise to cocycles 
by calculating $(wg)^3$. 
It then describes a full list of representatives of $H^1(\ZZ/3\ZZ,\cal{W})$ by simply calculating $r^{-1}(wg)r$ 
for every $r\in \cal{W}$ and every $w$ such that $(wg)^3 = \mathrm{id}$. 

It turns out that $H^1(\ZZ/3\ZZ,\cal{W})$ has exactly $2$
elements, with the non-trivial element being represented by taking
$$ 
w \,=\, \left(  \begin{matrix}
-1 & &  &  \\  -1&  & 1 & -1 \\ -1 & 1 &  & -1 \\-1  &  & 1& 
\end{matrix} \right)
$$
in matrix form. 
It therefore suffices to check that in this case $wg$ and $g$ are not conjugate in $\mathrm{GL}(\Lambda_F)$, 
which we can do by observing that $g$ has $1$ as an eigenvalue, but $wg$ does not, as they are all non-trivial cube roots of unity.
We therefore find that 
$$ 
H^1\left(\ZZ/3\ZZ,\cal{W}\right) \,\to\, H^1\left(\ZZ/3\ZZ,\mathrm{GL}(\Lambda_{F}) \right)
$$
has trivial kernel as claimed.
\qed\medskip

An important part of the induction step will be inflation-restriction in non-abelian cohomology. 
If we let $N\vartriangleleft G$ be a normal subgroup of $G$,
then we can consider the quotient $\mathbf{T}/N$, 
which is again a finite collection of finite Dynkin diagrams, on which $G/N$ acts. 
Moreover, we have, by Lemma \ref{lemma: GinvW}, isomorphisms
$$
\Lambda_{\mathbf{T}/N} \,\cong\, \Lambda^N_\mathbf{T},\;\; 
  \cal{W}_{\mathbf{T}/N} \,\cong\, \cal{W}^N_\mathbf{T},
$$
as well as a natural map
$$
\mathrm{GL}(\Lambda_{\mathbf{T},F})^N \,\to\, \mathrm{GL}(\Lambda^N_{\mathbf{T},F})
$$
which is $G/N$-equivariant.

\begin{Lemma}
  \label{lemma: infres} 
   If both maps
   $$
   \begin{array}{lcl}
      H^1\left(N,\cal{W}_\mathbf{T}\right) &\rightarrow&  H^1\left(N,\mathrm{GL}(\Lambda_{\mathbf{T},F})\right) \\
      H^1\left(G/N,\cal{W}_{\mathbf{T}/N}\right) &\rightarrow&  H^1\left(G/N,\mathrm{GL}(\Lambda_{\mathbf{T}/N,F})\right)
    \end{array}
   $$
have trivial kernel, then so does
$$
H^1\left(G,\cal{W}_\mathbf{T}\right) \,\to\, H^1\left(G,\mathrm{GL}(\Lambda_{\mathbf{T},F})\right).
   $$
\end{Lemma}

\prf 
By inflation-restriction \cite[Ch. I, \S5.8 (a)]{S94}, we have a commutative diagram of pointed sets
$$ 
\xymatrix{    
 H^1(G/N,\cal{W}_{\mathbf{T}/N}) \ar[r] \ar[d] &  H^1(G, \cal{W}_\mathbf{T}) \ar[r] \ar[d] &  H^1(N, \cal{W}_\mathbf{T}) \ar[d]  \\  
 H^1(G/N, \mathrm{GL}(\Lambda_{\mathbf{T},F})^N)\ar[d]  \ar[r] & H^1(G,\mathrm{GL}(\Lambda_{\mathbf{T},F})) \ar[r] &  H^1(N,\mathrm{GL}(\Lambda_{\mathbf{T},F}))   \\  
  H^1(G/N, \mathrm{GL}(\Lambda_{\mathbf{T}/N,F}))  }  
$$
with exact rows and injective left hand horizontal maps. 
The claim now follows from a simple diagram chase.
\qed\medskip

The other key component of the induction step is the following. Suppose $G$ acts on $\mathbf{T}=(T_1\ldots,T_d)$ in such a way that the action on the irreducible components $T_i$ is transitive. Let $H_i\leq G$ be the stabiliser of $T_i$, thus $H_i$ acts on $T_i$ and therefore on $\cal{W}_{T_i}$.

\begin{Lemma} 
 \label{lemma: induced} 
  For any $i\in \{1,\ldots,d\}$, there is an isomorphism
$$
\cal{W}_\mathbf{T} \,\cong\, \mathrm{Ind}_{H_i}^G\, \cal{W}_{T_i}
  $$
of $G$-modules.\end{Lemma}

\prf 
Let $\mathbf{T}'=(T_1,\ldots,\widehat{T}_i,\ldots,T_d)$, so that we have $\cal{W}_\mathbf{T}\cong \cal{W}_{T_i} \times \cal{W}_{\mathbf{T}'}$ as non-abelian $H_i$-modules. 
For all $1\leq j \leq d$ let $p_j:\cal{W}_{\mathbf{T}}\rightarrow \cal{W}_{T_j}$ denote the $j$th projection. Then, there is a natural map
$$
 \cal{W}_{\mathbf{T}} \,\to\, \mathrm{Ind}_{H_i}^G\, \cal{W}_{T_i}
$$
of $G$-modules defined by
\[ w\mapsto \left( f_w(x) := p_i(x(w)) \right). \]
These two $G$-modules have the same order, it therefore suffices to show that the map is injective. Since the action of $G$ on the irreducible components is transitive, $f_w(x)=1$ for all $x$ implies that $p_j(w)=1$ for all $j$ and thus $w=1$.
\qed\medskip

We can now complete the proof of Theorem \ref{theo: nagc main}.
\medskip

\prf[Proof of Theorem \ref{theo: nagc main}] 
We induct on the order of $G$, if $\lvert G\rvert =1$ then there is nothing to prove.
For the induction step we break up $\mathbf{T}=(\mathbf{T}_1,\ldots,\mathbf{T}_r)$ 
into orbits for the action of $G$ on the irreducible components $T_i$. 
Then, we have corresponding decompositions 
$$
\cal{W}_\mathbf{T} \,\cong\, \prod_{j=1}^r\cal{W}_{\mathbf{T}_j} \mbox{ \quad and \quad }
  \Lambda_\mathbf{T} \,\cong\, \bigoplus_{j=1}^r\Lambda_{\mathbf{T}_j} 
$$
and thus,
$$
H^1\left(G,\cal{W}_\mathbf{T}\right) \,\cong\,
  \prod_{j=1}^r H^1\left(G,\cal{W}_{\mathbf{T}_j}\right). 
$$
It therefore suffices to show that each 
$$
H^1\left(G,\cal{W}_{\mathbf{T}_j}\right) \,\to\, 
  H^1\left(G, \mathrm{GL}(\Lambda_{\mathbf{T},F})\right)
$$
has trivial kernel (note that at this point we really mean $\mathbf{T}$ on the RHS, not $\mathbf{T}_j$). 
The decomposition
$$
\Lambda_{\mathbf{T},F} \,\cong\, \bigoplus_j\, \Lambda_{\mathbf{T}_j,F}
$$
is invariant under both $G$ and $\cal{W}_{\mathbf{T}_j}$, and moreover, any element of 
$\cal{W}_{\mathbf{T}_j}$ acts trivially on $\bigoplus_{j'\neq j} \Lambda_{\mathbf{T}_{j'},F}$. 
By general theory, $\alpha$ mapping to the trivial class in 
$H^1(G, \mathrm{GL}(\Lambda_{\mathbf{T},F} ))$ is equivalent to the existence 
of an isomorphism
$$
\bigoplus_j \Lambda_{\mathbf{T}_j,F} \,\cong\, \bigoplus_j \Lambda^\alpha_{\mathbf{T}_j,F} 
$$
of $G$-representations, where on the RHS we have ``twisted'' the $G$-action by $\alpha$. 
Similarly, $\alpha$ mapping to the trivial class in $H^1(G, \mathrm{GL}(\Lambda_{\mathbf{T}_j,F} ))$ 
is equivalent to the existence of an isomorphism
$$
\Lambda_{\mathbf{T}_j,F} \,\cong\,  \Lambda^\alpha_{\mathbf{T}_j,F}.
$$
Since $\mathrm{char}(F)=0$, the category of $G$-representations in finite dimensional $F$-vector spaces is semi-simple, and as
$$
\Lambda_{\mathbf{T}_{j'},F} \cong  \Lambda^\alpha_{\mathbf{T}_{j'},F}
$$
for $j\neq j'$, we see that
$$
\Lambda_{\mathbf{T}_j,F} \,\cong\,  \Lambda^\alpha_{\mathbf{T}_j,F}
$$
if and only if
$$
\bigoplus_j \Lambda_{\mathbf{T}_j,F} \,\cong\, \bigoplus_j \Lambda^\alpha_{\mathbf{T}_j,F}.
$$
From this we deduce that $\alpha$ is trivial in $H^1(G, \mathrm{GL}(\Lambda_{\mathbf{T},F} ))$ 
if and only if it is so in $H^1(G, \mathrm{GL}(\Lambda_{\mathbf{T}_j,F} ))$. 
Thus, it suffices to show that each map
$$
H^1\left(G,\cal{W}_{\mathbf{T}_j}\right) \,\to\, H^1\left(G, \mathrm{GL}(\Lambda_{\mathbf{T}_j,F} )\right)
$$
has trivial kernel. 
Hence, we can assume that the action on the irreducible components of $\mathbf{T}=(T_1,\ldots,T_d)$ is transitive.

If $d>1$, then we use Lemma \ref{lemma: induced}. The stabilizer $H_i\leq G$ of some $T_i$ is a \emph{proper} subgroup of $G$, and by \cite[Ch. I, \S5.8 (b)]{S94} we know that 
$$ 
H^1\left(G,\cal{W}_{\mathbf{T}}\right) \,\isomto\, H^1\left(H_i,\cal{W}_{T_i}\right). 
$$
In particular, if we let $\mathbf{T}'=(T_1,\ldots,\widehat{T}_i,\ldots,T_d)$, we see that the pull-back map
$$
H^1\left(G,\cal{W}_{\mathbf{T}}\right) \,\to\, 
H^1\left(H_i,\cal{W}_{\mathbf{T}}\right) \,\isomto\, 
H^1\left(H_i,\cal{W}_{T_i}\right) \times H^1\left(H_i,\cal{W}_{\mathbf{T}'}\right) 
$$
has trivial kernel. By the induction hypothesis we know that
$$
H^1\left(H_i,\cal{W}_{\mathbf{T}}\right)  \,\to\, H^1\left(H_i, \mathrm{GL}(\Lambda_{\mathbf{T},F} )\right) 
$$
has trivial kernel, thus so does 
$$
H^1\left(G,\cal{W}_{\mathbf{T}}\right) \,\to\, H^1\left(G, \mathrm{GL}(\Lambda_{\mathbf{T},F} )\right) 
$$
by a simple diagram chase. 

If $d=1$, then $\mathbf{T}$ has a unique irreducible component. By inspecting
the kernel of the $G$-action and using Lemma \ref{lemma: infres}, we may consider separately the cases when the $G$-action is faithful and when it is 
trivial. If the $G$-action is trivial, we can identify $H^1(G,\cal{W}_\mathbf{T})$ 
with the set of conjugacy classes of homomorphisms 
$G\rightarrow \cal{W}_\mathbf{T}$,
and $H^1(G,\mathrm{GL}(\Lambda_{\mathbf{T},F}))$ with the set of conjugacy 
classes of homomorphisms 
$G\rightarrow \mathrm{GL}(\Lambda_{\mathbf{T},F})$. But since $\cal{W}_\mathbf{T}\rightarrow  \mathrm{GL}(\Lambda_{\mathbf{T},F})$ is injective, we may therefore conclude by observing that the only homomorphism $G\rightarrow \mathrm{GL}(\Lambda_{\mathbf{T},F})$ conjugate to the trivial homomorphism is the trivial homomorphism itself.

If the $G$-action is faithful (and $G$ is non-trivial), then $G$ is either $\ZZ/2\ZZ$, $\ZZ/3\ZZ$ or $S_3$. Once more applying Lemma \ref{lemma: infres} 
we can reduce either to Proposition \ref{prop: nagc faith}, or to the already handled case when the $G$-action is trivial. 
\qed\medskip

\section{Crystalline Galois representations and K3 surfaces}
\label{sec: crys rep}

In this section, we establish a $p$-adic criterion for potential good reduction
of K3 surfaces (Theorem \ref{thm: first main result in introduction}),
as well as the counter-examples (Remark \ref{rem: counter examples})
from the introduction.

\subsection{A \texorpdfstring{$\bm{p}$}{p}-adic criterion for potential good reduction of K3 surfaces}
We start with the proof of Theorem \ref{thm: first main result in introduction}.
\medskip

\prf
The equivalences $(1)\Leftrightarrow(2)\Leftrightarrow(3)$ have already been established in 
\cite[Theorem 1.3 and Corollary 1.4]{LM18}.
For every prime $\ell$ (including $\ell=p$), we set
$V_\ell:=\Het{2}(X_{\overline{K}},\QQ_\ell)$ and let $\rho_\ell:G_K\to {\rm GL}(V_\ell)$ 
be the $\ell$-adic Galois representation arising from geometry.
The $\QQ_\ell$-dimension of $V_\ell$ is independent of $\ell$ and
equal to the second Betti number $b_2:=b_2(X)$.

Next, assume that $X$ has good reduction after some finite and 
unramified extension $K\subseteq K'$.
By Theorem \ref{thm:crysrep}, the restriction of $\rho_p$ to $G_{K'}$ is 
a crystalline $G_{K'}$-representation.
But then, $\rho_p$ is a crystalline $G_K$-representation
by Lemma \ref{lemma: change of group} (\ref{num: lcog1}).
This establishes $(1)\Rightarrow(4)$.

Finally, assume that $\rho_p$ is crystalline. Since the conditions (3) and (4) are equivalent over $K$ and $\widehat{K}$, we may assume that $K$ is complete.
By Lemma \ref{lemma: change of group} (\ref{num: lcog3}), the induced
$I_K$-action on $\DDpst(V_p)$ is trivial.
In particular, for every $g\in I_K$, the trace of 
the $K_0^{\rm nr}$-linear map $g^*$ on $\DDpst(V_p)$
is equal to $b_2$.
But then, also for every prime $\ell\neq p$, the trace of the induced
$\QQ_\ell$-linear map $g^*$ on 
$V_\ell$ is equal to $b_2$ by \cite[Theorem 3.1]{Och99}. Note that as long as we only care about the action of the inertia group $I_K$, Ochiai's theorem only requires the residue field $k$ to be perfect. We can also apply this result to deduce $\ell$-independence of traces for $H^2$, since we already know $\ell$-independence of traces for $H^0$ and $H^4$ (trivially), and both $H^1$ and $H^3$ vanish.

Now, since $\rho_p$ is crystalline, it follows from \cite[Theorem 1.1]{Mat15} that there exists 
a finite and possibly ramified extension $L/K$ such that $X_L$ has good reduction 
over $\OO_L$.
Thus, for every $\ell\neq p$ the image of $I_K$ in ${\rm GL}(V_\ell)$ is finite.
Since $\QQ_\ell$ is of characteristic zero, a linear automorphism of 
finite order on $V_\ell$ is trivial if and only if its trace is equal to $b_2$.
This proves that $I_K$ acts trivially on $V_\ell$ for all $\ell\neq p$
and establishes $(4)\Rightarrow(3)$.
\qed\medskip

\begin{Remark}
 Let us comment on Theorem \ref{thm: first main result in introduction} and
 sketch some alternative strategies to prove (parts of it):
  \begin{enumerate} 
 \item It seems plausible that one can prove the equivalence $(1)\Leftrightarrow(4)$
  by adapting the $\ell$-adic arguments in the proof of
  \cite[Theorem 6.1]{LM18} to the $p$-adic situation.
 \item
  If $A$ is an Abelian variety over $K$, then 
  the analogous equivalences $(2)\Leftrightarrow(3)\Leftrightarrow(4)$ 
  for $\Het{1}(A_{\overline{K}},\QQ_\ell)$  
  follow from work of Serre and Tate \cite{ST68} in case $\ell\neq p$
  and from work of Fontaine \cite{Fon79}, 
  Mokrane \cite{Mok93}, and Coleman and Iovita \cite{CI99} 
  for $\ell=p$.

  It seems plausible that one might be able to use this result 
  together with the Kuga--Satake correspondence \cite{KS67, And96, Riz10} 
  to give another proof of the equivalences  $(2)\Leftrightarrow(3)\Leftrightarrow(4)$ for K3 surfaces.

 \item Quite generally, the equivalences
  $(2)\Leftrightarrow(3)\Leftrightarrow(4)$ are expected to hold
  for every smooth and proper variety over 
  $K$.
  \end{enumerate} 
\end{Remark}

\subsection{Counter-examples}

Next, as already noted in Remark \ref{rem: counter examples},
we show that the unramified extension stated 
in Theorem \ref{thm: first main result in introduction} (\ref{num :1st main in intro 1})
may be \emph{non-trivial}.
More precisely, the examples from \cite[Section 7]{LM18} 
provide the desired counter-examples.

\begin{Theorem}
 \label{thm: counter examples}
 For every prime $p\geq5$, there exists a K3 surface $X(p)$ over $\QQ_p$ such that
  \begin{enumerate} 
  \item the $G_{\QQ_p}$-representation on  $\Het{2}(X_{\overline{\QQ}_p},\QQ_\ell)$ is
    unramified for all $\ell\neq p$.
  \item the $G_{\QQ_p}$-representation on $\Het{2}(X(p)_{\overline{\QQ}_p},\QQ_p)$ is crystalline,
  \item $X(p)$ has good reduction over the unramified extension $\QQ_{p^2}$, but
  \item $X(p)$ does not have good reduction over $\QQ_p$.
  \end{enumerate} 
\end{Theorem}

\prf
Let $X(p)$ be one of the examples of \cite[Example 7.1]{LM18}.
By \cite[Theorem 7.2]{LM18}, it satisfies claims (1), (3) and (4).
Since it has good reduction over $\QQ_{p^2}$, the $G_{\QQ_{p^2}}$-representation on 
$V_p:=\Het{2}(X_{\overline{\QQ}_p},\QQ_p)$ is crystalline.
By Lemma \ref{lemma: change of group} (1), also the $G_{\QQ_{p}}$-representation on 
$V_p$ is crystalline, which establishes claim (2) for $X(p)$.
\qed\medskip

\subsection{Enriques surfaces}
By definition, an {\em Enriques surface} over a field $F$ 
is a smooth and proper surface $S$ over $F$ with $\omega_S\equiv\OO_S$ 
and $b_2(S)=10$, where $\equiv$ denotes numerical equivalence.
If ${\rm char}(F)\neq2$, then there exists a unique geometrically connected
\'etale cover $\widetilde{S}\to S$ of degree $2$.
In fact, $\widetilde{S}$ is a K3 surface over $F$ and it is called the
{\em K3 double cover} of $S$.

The following result is analogous to \cite[Lemma 3.4]{LM18}, which is why we leave
its proof to the reader.

\begin{Lemma}
 Assume that $p>0$, and let $Y$ be an Enriques surface over $K$.
 Then, there exists a finite extension $L/K$ such that the $G_L$-representations
 on $\Het{n}(X_{\overline{K}},\QQ_p)$ are crystalline for all $n$.
\end{Lemma}

Moreover, combining Theorem \ref{thm: first main result in introduction}
and \cite[Theorem 3.6]{LM18}, we obtain the following result.

\begin{Theorem}
  For every prime $p\geq5$, there exists an Enriques surface $Y(p)$ over 
  $\QQ_p$, such that
   \begin{enumerate} 
   \item the K3 double cover $X(p)$ of $Y(p)$ has good reduction,
   \item the $G_{\QQ_p}$-representation on $\Het{2}(X(p)_{\overline{\QQ}_p},\QQ_p)$
    is crystalline,
    \item $Y(p)$ has semi-stable reduction of flower pot type, but
    \item $Y(p)$ does not have potential good reduction.
   \end{enumerate} 
\end{Theorem}

In particular, there is no analog of Theorem \ref{thm: first main result in introduction}
for Enriques surfaces, not even in terms of their K3 double covers.

\section{Canonical reductions and RDP models of pairs} 
\label{sec: rdpcr}

If $X$ is a K3 surface over $K$, satisfying $(\star)$ and the equivalent conditions of Theorem \ref{thm: first main result in introduction}, then we may appeal to \cite[Theorem 1.3]{LM18} to deduce the existence of ``reasonably nice'' models of $X$ over $\OO_K$.

\begin{Theorem}[Liedtke, Matsumoto]
  \label{thm: RDP models} 
  Let $X$ be a K3 surface over $K$ that satisfies $(\star)$ and
  the equivalent conditions of Theorem \ref{thm: first main result in introduction}.
  Then, there exists a projective and flat model
  $$ 
  \cal{X} \,\rightarrow\, \Spec \OO_K
  $$
 for $X$ such that:
   \begin{enumerate} 
    \item the special fiber $\cal{X}_k$ has at worst canonical singularities, and
    \item its minimal resolution is a K3 surface over $k$.
   \end{enumerate} 
\end{Theorem}

While the ``RDP model'' provided by this theorem is not unique, 
it is possible to construct, once we fix a polarization $\cal{L}$ on $X$, 
a completely canonical such model $P(X,\cal{L})$ as follows: 
first of all, one takes a finite, unramified extension $L/K$, Galois with group $G$, 
such that $X$ has good reduction over $L$. 
Let $k_L$ denote the residue field of $L$. 
Then, for any smooth model $\cal{Y}\rightarrow \Spec\OO_L$ for $X_L$ 
we can find a birational map $\cal{Y}\dashrightarrow \cal{Y}^+$ to another smooth model $\mathcal{Y}^+$ for $X_L$, 
such that the specialization $\cal{L}_{k_L}^+$ of $\cal{L}_L$ 
on the special fiber $\cal{Y}_{k_L}^+$ is big and nef, \cite[Proposition 4.5]{LM18}.
The projective scheme
$$ 
P(X_L,\cal{L}_L) \,:=\, \mathrm{Proj}\left(  \bigoplus_{n\geq0} H^0(\cal{Y}^+,\cal{L}^{+,\otimes n}) \right)
$$
only depends on $X_L$ and $\cal{L}_L$ \emph{up to a canonical isomorphism} 
by the proof of \cite[Proposition 4.7]{LM18}. 
Moreover, the rational $G$-action on $X_L$ extends to a regular $G$-action on 
$P(X_L,\cal{L}_L)$ by \cite[Proposition 5.1]{LM18}, and we define
$$ 
P(X,\cal{L}) \,=\, P(X_L,\cal{L}_L)/G
$$
to be the quotient. 
This is then a flat, projective $\OO_K$-scheme whose special fibre $P(X,\cal{L})_k$ 
has canonical singularities \cite[Proposition 4.6]{LM18}. 
The model $P(X,\cal{L})$ depends only on $X$ and $\cal{L}$, again up to a \emph{canonical} isomorphism. 
Moreover, the minimal resolution of singularities
$$ 
Y \,\rightarrow\, P(X,\cal{L})_k
$$
of the special fibre is a K3 surface over $k$ by \cite[Proposition 4.6]{LM18}. This K3 surface $Y$ is unique up to canonical isomorphism and depends only on $X$
by a theorem of Matsusaka and Mumford \cite{MM64}
and the theory of minimal
models, see also \cite[Proposition 4.7]{LM18}.

\begin{Definition} 
 We call $Y$ the \emph{canonical reduction} of $X$
 and $P(X,\cal{L})$ the \emph{(canonical) RDP model} of the pair $(X,\cal{L})$.
\end{Definition}

We note that all of this is compatible with base change, in that we have
$$
P(X_L,\cal{L}_L) \,\cong\, P(X,\cal{L}) \otimes_{\OO_K} \OO_L
$$
for \emph{any} finite extension $L/K$. 
This then implies that the canonical reduction of $X_L$ is just the base change 
of the canonical reduction of $X$. 
Moreover, when $X$ \emph{does} have good reduction, 
then the canonical reduction $Y$ is the special fibre of 
\emph{any} smooth model of $X$, again by the results of \cite{MM64}.

\section{The Weyl group of a K3 surface with potentially good reduction}\label{sec: weyl k3}

Let $X/K$ be a K3 surface, satisfying $(\star)$ and the equivalent conditions of 
Theorem \ref{thm: first main result in introduction}. 
Thus, $X$ admits good reduction over a finite and unramified Galois 
extension $L/K$, that is, there exists a smooth and proper model 
$\cal{Y}\rightarrow \cal{O}_L$ for $X_L$. 
If we let $G$ denote the Galois group of $L/K$, then there is a 
natural semi-linear $G$-action on $X_L$, and hence a \emph{rational} $G$-action 
on $\cal{Y}$. 
The following is an analogue of \cite[Corollary 5.12]{LM18} in the unramified case.

\begin{Proposition} 
 \label{prop: quotient} 
 Assume that the rational and semi-linear $G$-action on $\cal{Y}$ is in fact regular. 
 Then, the quotient $\cal{X}:=\cal{Y}/G$ exists as a smooth and proper algebraic space 
 over $K$, and we have a $G$-equivariant isomorphism 
 $\cal{Y}\cong \cal{X}\otimes_{\cal{O}_K} \cal{O}_L$.
\end{Proposition}

\prf 
The $G$-action is equivalent to descent data for $\cal{Y}$ via the finite \'etale cover 
$\Spec \cal{O}_L\rightarrow \Spec \cal{O}_K$, and for algebraic spaces 
all such descent data are effective. 
Moreover, smoothness and properness can both be checked after taking such a cover.
\qed\medskip

In particular, $X$ has good reduction over $K$ itself if and only if such a $\cal{Y}$ 
can be chosen for which the rational $G$-action is in fact regular. 
Thus, a natural question arises of whether or not we can give conditions 
under which this will happen. 
This will involve a detailed study of the action of ``flops'' of arithmetic 3-folds on 
the Picard group of the special fibre of $\cal{Y}$, similar to that carried out in \cite{LM18} for $\ell$-adic cohomology. 
We will then use this to provide a cohomological obstruction for the regularity 
of the $G$-action (Proposition \ref{prop: NOS1} below).

\subsection{Flops and the Weyl group}\label{subsec: flops weyl}

We start by examining the action of flops on smooth models of a 
fixed, polarized K3 surface $(X,\cal{L})$ over $K$, again assuming $(\star)$ 
and the conditions of Theorem \ref{thm: first main result in introduction}. 
In this case, we have the canonical RDP model $P(X,\cal{L})$ over $\OO_K$, 
as well as the canonical reduction $Y$, which is the minimal resolution of singularities
$$ 
 Y \,\to\,  P(X,\cal{L})_k
$$
of the special fiber of $P(X,\cal{L})$. 
We will denote by
\begin{align*} 
E_{X,\cal{L}} &\subset Y \\
\Lambda_{X,\cal{L}} & \subset \mathrm{Pic}(Y) \\
\cal{W}_{X,\cal{L}} &\leq \mathrm{Aut}_{\ZZ}(\mathrm{Pic}(Y))
\end{align*}
respectively the exceptional locus of $Y \rightarrow P(X,\cal{L})_k$, and the associated 
root lattice and Weyl group as in Section \ref{sec: root}. 
We will call $\Lambda_{X,\cal{L}}$ and $\cal{W}_{X,\cal{L}}$ the \emph{root lattice}
and the \emph{Weyl group} of the pair $(X,\cal{L})$.

\begin{Definition} 
 We say that a smooth model $\cal{X}$ of $X$ is $\cal{L}$-terminal, or sometimes that 
 $(\cal{X},\cal{L})$ is terminal, or even that $\cal{X}$ is a terminal model of $(X,\cal{L})$, 
 if the specialization $\cal{L}_k$ on the special fibre $\cal{X}_k$ is big and nef.
\end{Definition}

If we have a birational map $f:\cal{X}\dashrightarrow \cal{X}^+$ of smooth models 
of K3 surfaces, with $(\cal{X},\cal{L})$ terminal, then we can consider the push-forward  $\cal{L}^+$ of $\cal{L}$ on the generic fiber of $\cal{X}^+$. 

\begin{Definition} 
We say that $f$ is a \emph{terminal} birational map if $(\cal{X}^+,\cal{L}^+)$ is also terminal.
\end{Definition}

If $X$ has good reduction over $K$, that is, if \emph{some} smooth model $\cal{X}$ exists, 
then an $\cal{L}$-terminal model exists by \cite[Proposition 4.5]{LM18}. 
We note that if $\cal{X}$ is an $\cal{L}$-terminal model for $X$, then by \cite[Proposition 4.6(2)]{LM18} 
we can identify $E_{X,\cal{L}}$ with exactly those curves $E\subset \cal{X}_k$ 
on the special fibre that satisfy $\cal{L}_k\cdot E=0$.

Given a model $\cal{X}$ of $X$, and a birational map 
$f:\cal{X}\dashrightarrow \cal{X}^+$ to some other smooth model of a K3 surface 
$X^+/K$, we may consider the (closure of the) graph of $f$
$$
\Gamma_f \,\subset\, \cal{X}\times_{\OO_K} \cal{X}^+,
$$
as well as its special fibre $\Gamma_{f,k}\subset Y\times_k \cal{X}^+_k$. 
Note that the generic fiber $\Gamma_{f_K}$ is simply the graph of the induced isomorphism 
$f_K:X\isomto X^+$, but the same is not necessarily true for $\Gamma_{f,k}$. 
This latter cycle induces a homomorphism
$$
\begin{array}{ccccc}
   \tilde{s}_f &:&  \mathrm{Pic}(\cal{X}_k^+) &\rightarrow& \mathrm{Pic}(Y)   \\
     &&D &\mapsto& p_{1*}\left(\Gamma_{f,k} \cap p_2^*D\right) 
 \end{array}
$$
on Picard groups, where $p_i: Y\times_k \cal{X}_k^+ \rightrightarrows Y,\cal{X}_k^+$ 
are the two projections. 
Composing with the push-forward via the induced isomorphism $f_k:Y\rightarrow \cal{X}_k^+$,
we obtain an endomorphism
$$ 
s_f \,:\,  \mathrm{Pic}(Y) \,\to\, \mathrm{Pic}(Y),
$$
which can be described as the map induced by the pull-back cycle 
$\Gamma\subset Y\times_k Y$ via the same formula.

\begin{Proposition}  
 \label{prop: sf first props}  
  \begin{enumerate} 
   \item \label{sffp1} The map $s_f$ preserves the intersection pairing, that is,
$$
s_f(D_1) \cdot s_f(D_2) \,=\, D_1\cdot D_2
 $$
for $D_i\in \mathrm{Pic}(Y)$. 
\item \label{sffp2}  If $f,g$ are composable birational maps, then we have 
$$
s_{g\circ f } \,=\, s_f \circ f_k^* \circ s_g \circ (f_k^*)^{-1}.
  $$
\item \label{sffp3} The map $s_f$ is invertible, and we have
$$
     (s_f)^{-1} \,=\, f_k^* \circ s_{f^{-1}} \circ (f_k^*)^{-1}.
$$
\end{enumerate} 
\end{Proposition}

\prf 
Since clearly $s_\mathrm{id} = \mathrm{id}$, we see that (\ref{sffp2})$\Rightarrow$(\ref{sffp3}), 
and to prove (\ref{sffp1}) and (\ref{sffp2}) it suffices to work instead with 
$\tilde{s}_f$, that is, to show that $\tilde{s}_f$ preserves the intersection pairing, 
and that $\tilde{s}_{g\circ f} = \tilde{s}_f\circ \tilde{s}_g$. 

To prove this, we will anticipate somewhat a later argument and use cohomology to give another interpretation of $\tilde{s}_f$. Fix a prime $\ell\neq p$, and extend $\tilde{s}_f$ to an endomorphism on $\ell$-adic cohomology
$$
\begin{array}{ccccc}
  \tilde{s}_{f,\ell} &:&  \Het{2}(\cal{X}^+_{\bar k},\QQ_\ell) &\rightarrow& \Het{2}(Y_{\bar k},\QQ_\ell)  \\
    &&\alpha &\mapsto& p_{1*}\left([\Gamma_{f,k}] \cup p_2^*\alpha\right) 
 \end{array}
$$
using the same formula. 
It therefore suffices to show that $\tilde{s}_{f,\ell}$ is compatible with composition and preserves 
the Poincar\'e pairing. 
This can be seen by applying \cite[Lemma 5.6]{LM18}, which implies that there is a 
commutative diagram
$$ 
\xymatrix{  
  \Het{2}(X^+_{\overline{K}},\QQ_\ell)  \ar[r]^-{\cong}\ar[d]_{f_K^*} & \Het{2}(\cal{X}^+_{\bar k},\QQ_\ell) \ar[d]^{\tilde{s}_{f,\ell} } \\ 
  \Het{2}(X_{\overline{K}},\QQ_\ell) \ar[r]^-{\cong} & \Het{2}(Y_{\bar k},\QQ_\ell)  
} 
$$
where the horizontal arrows come from smooth and proper base change, 
and the left hand vertical arrow is simply pull-back along the map 
$f_K:X \rightarrow X^+$ on generic fibers.
\qed\medskip

The main result we want to prove is that when $f$ is terminal, $s_f$ is an element of the Weyl group
$$
\cal{W}_{X,\cal{L}} \,\leq\, \mathrm{Aut}_{\ZZ}\left(\mathrm{Pic}(Y)\right).
$$
First we will need to show that we can \emph{construct} elements of the Weyl group 
via this method. 

\begin{Proposition} 
 \label{prop: surjW} 
 Let $(\cal{X},\cal{L})$ be terminal and assume that all irreducible components of 
 $E_{X,\cal{L}}$ are geometrically irreducible. 
 Then, for any element $w \in \cal{W}_{X,\cal{L}}$, there exists a terminal birational map
$$
f \,:\, \cal{X} \,\dashrightarrow\, \cal{X}^+
 $$
such that $s_f=w$ as automorphisms of $\mathrm{Pic}(Y)$. 
\end{Proposition}

\prf 
By Proposition \ref{prop: sf first props}(\ref{sffp2}), it suffices to treat the case when 
$w=s_E$ is the reflection in some irreducible component $E\subset E_{X,\cal{L}}$. 
In this case, we follow the proof of \cite[Proposition 4.2]{LM18}. 

First of all, there exists a proper and birational morphism 
$\cal{X}\rightarrow \cal{X}'$ that contracts $E$ and nothing else. 
Let $y\in \cal{X}'$ denote the image of $E$, $\widehat{\cal{X}}'_y=\mathrm{Spf}(R)$ the formal completion at $y$, and $\widehat{\cal{X}}\rightarrow \widehat{\cal{X}}'_y$ the formal fiber. 
Again, following \cite[Proposition 4.2]{LM18}, we can find an automorphism 
$t:\widehat{\cal{X}}'_y\rightarrow \widehat{\cal{X}}'_y$ which induces $-1$ 
on the Picard group $\mathrm{Pic}(R)$. 
Then, letting $\widehat{\cal{X}}^+\rightarrow \widehat{\cal{X}}'_y$ 
 denote the composite 
$\widehat{\cal{X}}\rightarrow \widehat{\cal{X}}'_y \overset{t}{\rightarrow} \widehat{\cal{X}}'_y$ we can algebraize to obtain a proper, birational morphism 
$\cal{X}^+\rightarrow \cal{X}'$ and hence a birational map 
$f:\cal{X}\dashrightarrow \cal{X}^+$. We claim that $s_f=s_E$.

Indeed, since $f$ induces an isomorphism
$$
  \cal{X} \setminus E \,\isomto\, \cal{X}^+ \setminus f_k(E),
$$
we note that the special fiber of the graph of $f$, considered as a cycle on 
$Y \times_k \cal{X}^+_k$, can be written as 
$$ 
\Gamma_{f,k} \,:=\, \Gamma_{f_k} + b(E\times f_k(E))
$$
for some $b\in \ZZ_{\geq 0}$. Thus, pulling back to $Y$ 
we can see that $s_f$ is given by the formula
$$ 
s_f(D) = D+b(E\cdot D) E. 
$$
By Proposition \ref{prop: sf first props}(\ref{sffp1}), we must have that $b\in \{0,1\}$. 
If $b=0$, then $f$ is an isomorphism and if $b=1$, then $s_f=s_E$.
It therefore suffices to show that
$$ 
  f\,:\,\cal{X}\,\dashrightarrow\, \cal{X}^+
$$
is \emph{not} biregular, and that $(\cal{X}^+,\cal{L}^+)$ is terminal. 
Since all these constructions are compatible with base change, 
to prove this we may therefore assume that $k=\bar k$ is algebraically closed (this is only to avoid any issues concerning the difference between geometric and closed points during the proof).

To see that $f$ is not biregular, it suffices to show that the dashed arrow in the commutative diagram
$$ \xymatrix{ 
  \widehat{\cal{X}} \ar[r] & \mathrm{Spf}(R) \\ 
  \widehat{\cal{X}} \ar[r] \ar@{-->}[u] & \mathrm{Spf}(R) \ar[u]_t }
$$
is not regular. 
Indeed, if it were, then we would have a commutative diagram
$$ \xymatrix{ 
  (\mathbf{R}^1g_*\cal{O}^*_{\widehat{\cal{X}}})_y \ar[d] \ar@{->>}[r] & \mathrm{Pic}(R) \ar[d]^{t^*} \\ 
  (\mathbf{R}^1g_*\cal{O}^*_{\widehat{\cal{X}}})_y \ar@{->>}[r] & \mathrm{Pic}(R)  }
$$
with $(\mathbf{R}^1g_*\cal{O}^*_{\widehat{\cal{X}}})_y$ generated by the class of 
the exceptional curve $E$.
In  particular, the left hand map would have to be the identity. 
This contradicts $t$ inducing $-1$ on $\mathrm{Pic}(R)$.

Finally, to see that $(\cal{X}^+,\cal{L}^+)$ is terminal, we simply note that since 
$\cal{L}_k\cdot E=0$, we can find a divisor representing $\cal{L}$ which doesn't meet $E$. 
Therefore, the specialization $\cal{L}^+_k$ is simply the push-forward 
$f_{k*}(\cal{L}_k)$ of $\cal{L}_k$ via the induced isomorphism $f_k:Y\rightarrow \cal{X}_k^+$.
\qed\medskip

With this in place, we can now prove that $s_f$ is an element of the Weyl group. 

\begin{Theorem} 
 \label{theo: flopref} 
 Let $\cal{X}$ be a terminal model of a polarized K3 surface $(X,\cal{L})$ over $K$
 and $f:\cal{X}\dashrightarrow \cal{X}^+$ a terminal birational map.
  \begin{enumerate} 
  \item \label{fr(1)} $s_f\in \cal{W}_{X,\cal{L}}\leq \mathrm{Aut}_{\ZZ}\left(\mathrm{Pic}(Y)\right)$.
  \item \label{fr(2)} $f$ is an isomorphism if and only if $s_f=\mathrm{id}$.
  \end{enumerate}  
\end{Theorem}

\prf 
If $L/K$ is a finite and unramified Galois extension with Galois group $G$, then 
$\cal{W}_{X,\cal{L}} \cong \cal{W}_{X_L,\cal{L}_L}^G$ by Corollary \ref{cor: GinvW}, 
and $f$ being an isomorphism can be detected over $L$. 
Hence, after possibly enlarging $K$, we can assume that all irreducible components $E_i\subset E_{X,\cal{L}}$ 
are geometrically irreducible. 
In particular, all singularities of $P(X,\cal{L})_k$ are rational double points of type ADE.

As in the proof of Proposition \ref{prop: surjW} above, we know that since $f$ is terminal, 
it induces an isomorphism
$$
\cal{X}\setminus E_{X,\cal{L}} \,\isomto\, \cal{X}^+ \setminus E_{X^+,\cal{L}^+},
$$
thus, we can describe the cycle $\Gamma\subset Y\times_k Y$ defining $s_f$ as
$$
\Gamma \,:=\, \Delta_{Y} \,+\, \sum_{i,j} b_{ij} (E_i \times E_j) 
$$
for some $b_{ij} \in \ZZ_{\geq0}$. Let $B=(b_{ij})\in \mathrm{M}_n(\ZZ)$ 
denote the corresponding matrix. 
We therefore have
$$
s_f (D) \,= \,D \,+\, \sum_{i,j} b_{ij} (D\cdot E_j)\cdot E_i
 $$
from which we see that $s_f$ preserves the root lattice
$$
\Lambda_{X,\cal{L}} \,\subseteq\, \mathrm{Pic}(Y).
$$
Moreover, since $f$ induces an isomorphism $P(X,\cal{L})\isomto P(X^+,\cal{L}^+)$, we deduce that 
$b_{ij}\neq 0$ only if $E_i$ and $E_j$ are in the same connected component of $E_{X,\cal{L}}$. 
Thus, $ s_f$ preserves the sub-lattice $\Lambda_l \subset \Lambda_{X,\cal{L}}$ 
corresponding to each connected component of $E_{X,\cal{L}}$. 
Using the formula $(s_f)^{-1} = f_k^* \circ s_{f^{-1}} \circ (f_k^*)^{-1}$,
 this is also the case for $(s_f)^{-1}$.

Hence, applying Proposition \ref{prop: sf first props}(\ref{sffp1}) we can write
$$
  s_f|_{\Lambda_{X,\cal{L}}} \,=\,\prod_l \tilde{w}_l
$$
as a product of elements $\tilde{w}_l \in \mathrm{O}(\Lambda_l)$ in the orthogonal group of each individual $\Lambda_l$. 
Moreover, by \cite[Theorem 1]{KM83}, we can write $\tilde{w}_l=\alpha_lw_l$ with $w_l$ in the Weyl group of 
$\Lambda_l$, and $\alpha_l$ coming from an automorphism of the Dynkin diagram of $\Lambda_l$. 
Now, let $w=\prod_l w_l\in \cal{W}_{X,\cal{L}}$. 
By Proposition \ref{prop: surjW}, there exists some birational map
$$ 
g\,:\, \cal{X}^+ \,\dashrightarrow\, \cal{X}^{+2}
$$
with $(\cal{X}^{+2},\cal{L}^{+2})$ terminal such that $s_g=(f_k^*)^{-1}\circ w^{-1}\circ  f_k^*$. 
Replacing $f$ by $g\circ f$ and applying Proposition \ref{prop: sf first props}(\ref{sffp2}), we may therefore assume that each $w_l=\mathrm{id}$, that is, that
$$ 
s_f|_{\Lambda_{X,\cal{L}}} \,=\,\prod_l \alpha_l
$$
is a product of elements coming from automorphisms of the relevant Dynkin diagrams.
We need to show that $s_f=\mathrm{id}$.

To see this, we note that the matrix of $s_f|_{\Lambda_{X,\cal{L}}}$ on a suitable basis of $\Lambda_{X,\cal{L}}$ is $I+BA$, where 
$A$ is the intersection matrix of $E_{X,\cal{L}}$ and $B$ is as above. 
Since both $s_f$ and $(s_f)^{-1}$ preserve each summand $\Lambda_l$, it follows that $B$ 
has to be a block sum matrix over the connected components of $E_{X,\cal{L}}$, the same is also true for $A$. 
Since we know that 
\[ (s_f|_{\Lambda_{X,\mathcal{L}}} -I)A^{-1} =B \in \mathrm{M}_n(\ZZ),\]
and everything appearing here decomposes as a relevant block sum over the connected components of $E_{X,\mathcal{L}}$ we may apply Lemma \ref{lemma: dynkaut} to conclude that in fact each $s_f|_{\Lambda_{X,\cal{L}}}=\mathrm{id}$. 
This then implies that $BA=0$, and since $A$ is negative definite, it therefore follows that $B=0$. 
Thus, we have $\Gamma=\Delta_Y$ and so $s_f=\mathrm{id}$.

The proof of (\ref{fr(2)}) is then entirely similar, since $s_f=\mathrm{id}$ implies that $BA=0$, whereas $f$ being regular is equivalent to having $B=0$. Negative definiteness of $A$ implies that $B=0\Leftrightarrow AB=0$.
\qed\medskip

This gives us a way of characterising terminal models of a given polarized K3 surface $(X,\cal{L})$ over $K$.

\begin{Corollary} 
 \label{cor: modweyl}  
 Let $(X,\cal{L})$ be a polarized K3 surface and assume that $X$ has good reduction over $K$. Then, the set 
 $$ 
\left\{ \cal{L}\text{-}\mathrm{terminal\;models\;of\;}X \right\}/\sim
 $$
of $\cal{L}$-terminal models up to equivalence is a torsor under the Weyl group $\cal{W}_{X, \cal{L}}$.
\end{Corollary}

\begin{Remark} Here, two such models $\cal{X}_1,\cal{X}_2$ of $X$ are said to be \emph{equivalent}, 
 if there exists an isomorphism $\cal{X}_1\isomto \cal{X}_2$ that induces the identity on the generic fiber $X$. 
 Put differently, this says that the rational map $\cal{X}_1\dashrightarrow \cal{X}_2$ defined by the identity 
 on the generic fiber is in fact regular. 
\end{Remark}

\prf 
Fix an $\cal{L}$-terminal model $\cal{X}$ for $X$. 
Then, we obtain a function
$$
\left\{ \cal{L}\text{-}\mathrm{terminal\;models\;of\;}X \right\}/\sim \,\to\, \cal{W}_{X,\cal{L}}
$$
as follows: given another $\cal{L}$-terminal model $\cal{X}^+$, we obtain a rational map
$$ 
f_{\cal{X}^+} \,:\, \cal{X} \,\dashrightarrow\, \cal{X}^+
$$
by taking the identity on generic fibers, and hence an element
$$
s_{\cal{X}^+} \,:=\,s_{f_{\cal{X}^+}}
$$
of the Weyl group. 
Injectivity of this map follows from Theorem \ref{theo: flopref}(\ref{fr(2)}) above, we would like to show surjectivity.

To do so, we will choose a finite and unramified Galois extension $L/K$, with Galois group $G$, such that all irreducible 
components of $E_{X_L,\cal{L}_L}$ are geometrically irreducible. 
If we have $w\in \cal{W}_{X,\cal{L}}$,  then applying Proposition \ref{prop: surjW} we know that there exists an 
$\cal{L}$-terminal model $\cal{Y}\rightarrow \OO_L$ of $X_L$ and a birational map 
$\cal{X} \otimes_{\OO_K} \OO_L \dashrightarrow \cal{Y}$, such that $s_{\cal{Y}}=w \in \cal{W}_{X_L,\cal{L}_L}$. 
But now, since $w$ is fixed by $G$, it follows from the already proved injectivity of 
$$
\left\{ \cal{L}_L\text{-}\mathrm{terminal\;models\;of\;}X_L \right\}/\sim \,\to\, \cal{W}_{X_L,\cal{L}_L}
$$
that so is $\cal{Y}$, for the natural $G$ action on the left hand side. 
Concretely, this means that the rational $G$-action on $\cal{Y}$ is regular, and thus, by Proposition \ref{prop: quotient},
we know that $\cal{Y}$ descends to $\OO_K$. 
\qed\medskip

Since the Weyl group is finite, we obtain the following rather modest corollary.

\begin{Corollary} 
 Any polarized K3 surface over $K$ with good reduction has finitely many 
 terminal models over $\cal{O}_K$.
\end{Corollary}

\subsection{Obstructions to good reduction} \label{subsec: ogr}

We will now use these results to study the relationship between good reduction and potential good reduction of a K3 surface 
over $K$. 
We will therefore take a polarized K3 surface $(X,\cal{L})$ over $K$, as always satisfying $(\star)$, 
such that the equivalent conditions of Theorem \ref{thm: first main result in introduction} hold. 

We thus know that $X$ admits good reduction over a finite and unramified Galois extension $L/K$, say with Galois group $G$ 
and residue field extension $k_L/k$. 
Let $P(X,\cal{L})$ denote the canonical RDP model over $\OO_K$ and $Y$ the canonical reduction of $X$ 
as in \S\ref{sec: rdpcr}. 
We therefore have 
$$
  \cal{W}_{X_L,\cal{L}_L}\, \leq\, \mathrm{Aut}_{\ZZ}\left(\mathrm{Pic}(Y_{k_L}) \right),
$$
which is invariant under the natural $G$-action on the latter. For any object over $L$ or $k_L$, and any $\sigma\in G$, we will denote by $(-)^\sigma$ the base change by $\sigma$, thus we have $(-)^{\sigma\tau} = ((-)^\tau)^\sigma$.

For any $\cal{L}_L$-terminal model $\cal{Y}$ of $X_L$ over $\OO_L$, 
we obtain a rational and semilinear $G$-action on the 
\emph{pair} $(\cal{Y},\cal{L}_L)$. 
Hence, by the results of Section \ref{subsec: flops weyl} we can define a function
$$ 
\alpha_\cal{Y} \,:\, G\,\to\, \cal{W}_{X_L,\cal{L}_L}
$$
as follows: for any $\sigma\in G$ we base change $\cal{Y}$ by $\sigma$ to obtain $\cal{Y}^\sigma$.
The regular $G$-action on the generic fibre $X_L$ provides a terminal (with respect to $\cal{L}_L$) birational map
$$
f_\sigma \,:\, \cal{Y} \,\dashrightarrow\, \cal{Y}^\sigma,
$$
which is $\OO_L$-linear. Then the induced rational map
\[ f_{\sigma,k_L} \,:\, Y_{k_L} \,\rightarrow\, Y_{k_L}^\sigma \]
on the special fiber is an isomorphism, and we can concretely describe the $G$-action on $\cal{W}_{X_L,\cal{L}_L}$ by the formula
\[ \sigma(s) = f_{\sigma,k_L}^*  \circ s^\sigma \circ (f_{\sigma,k_L}^*)^{-1}.\]
By Theorem \ref{theo: flopref} we may define an element $s_{f_\sigma}\in \cal{W}_{X_L,\cal{L}_L}$ of the Weyl group associated to $f_\sigma$, and we define a map
$$
\begin{array}{ccccc}
   \alpha_\cal{Y} &:& G &\rightarrow& \cal{W}_{X_L,\cal{L}_L}  \\
   &&\sigma &\mapsto& s_{f_\sigma}.
 \end{array}
$$
Then, we have the following crucial observation.

\begin{Proposition} 
 \label{prop: cocyc}  
 The map
 $$ 
   \alpha_\cal{Y} \,:\, G \,\to\, \cal{W}_{X_L,\cal{L}_L}
 $$ 
  is a 1-cocycle for the $G$-action on $\cal{W}_{X_L,\cal{L}_L}$. 
  Moreover, the $G$-action on $\cal{Y}$ is regular if and only if this cocycle is trivial, that is,
  satisfies $\alpha_\cal{Y}(\sigma)=1$ for all $\sigma\in G$. 
\end{Proposition}

\prf 
For the first claim, we use the fact that $f_{\sigma\tau} =f_\tau^\sigma \circ f_{\sigma}$, and simply calculate
\begin{align*} 
 \alpha_\cal{Y}(\sigma\tau) &= s_{f_{\sigma\tau}} =  s_{f_\tau^\sigma \circ   f_\sigma}  =s_{f_\sigma}\circ  f_{\sigma,k_L}^* \circ s_{f^\sigma_\tau}  \circ (f_{\sigma,k_L}^*)^{-1}  \\
  &=  s_{f_\sigma} \circ  f_{\sigma,k_L}^* \circ  s^\sigma_{f_\tau} \circ  (f_{\sigma,k_L}^*)^{-1} = s_{f_\sigma} \circ  \sigma(s_{f_\tau})  \\
  &= \alpha_\cal{Y}(\sigma)\circ  \sigma(\alpha_\cal{Y}(\tau))
\end{align*}
using Proposition \ref{prop: sf first props}(\ref{sffp2}). 
The second claim is an immediate consequence of Theorem \ref{theo: flopref}.
\qed\medskip

The next obvious question to ask is how the cocycle $\alpha_\cal{Y}$ changes as the $\cal{L}_L$-terminal model 
$\cal{Y}$ changes. 
So, let us suppose that we have another $\cal{L}_L$-terminal model $\cal{Y}^+$ for $X_L$, 
and let $f:\cal{Y}\dashrightarrow \cal{Y}^+$ be the rational map given by the identity on the generic fibers. 
Let $\alpha_{\cal{Y}^+}:G\rightarrow \cal{W}_{X_L,\cal{L}_L}$ denote the cocycle associated to the 
rational $G$-action on $\cal{Y}^+$, and let $s_f\in \cal{W}_{X_L,\cal{L}_L}$ 
be the element of the Weyl group associated to $f$, as provided by Theorem \ref{theo: flopref}.

\begin{Proposition} 
 \label{prop: cobound} 
  For all $\sigma\in G$ we have
$$
\alpha_{\cal{Y}^+}(\sigma) \,=\, s_f^{-1} \circ \alpha_{\cal{Y}}(\sigma) \circ \sigma(s_f).
  $$
\end{Proposition}

\prf 
If we view $\alpha_\cal{Y}$ as a function
\[\alpha_\cal{Y}: G\rightarrow \mathrm{Aut}_{\ZZ}(\mathrm{Pic}(\cal{Y}_{k_L})), \]
and $\alpha_{\cal{Y}^+}$ as a function
\[\alpha_{\cal{Y}^+}: G\rightarrow \mathrm{Aut}_{\ZZ}(\mathrm{Pic}(\cal{Y}^+_{k_L})), \]
then what we need to show is that 
$$
f_{k_L} ^* \circ \alpha_{\cal{Y}^+}(\sigma) \circ (f_{k_L}^*)^{-1} \,=\, s_f^{-1} \circ \alpha_{\cal{Y}}(\sigma) \circ \sigma(s_f).
  $$
Since the diagram
\[\xymatrix{ \cal{Y} \ar@{-->}[r]^-{f_\sigma} \ar@{-->}[d]_f & \cal{Y}^\sigma \ar@{-->}[d]^{f^\sigma} \\ \cal{Y}^+ \ar@{-->}[r]^-{f^{+}_\sigma} & \cal{Y}^{+,\sigma}}  \]
commutes for all $\sigma\in G$, we may again use Proposition \ref{prop: sf first props}(\ref{sffp2}) to calculate
\begin{align*}
\alpha_\cal{Y}(\sigma)\circ \sigma(s_f) &= s_{f_\sigma} \circ f_{\sigma,k_L}^* \circ s_{f}^\sigma \circ (f_{\sigma,k_L}^*)^{-1} \\
&= s_{f_\sigma} \circ f_{\sigma,k_L}^* \circ s_{f^\sigma} \circ (f_{\sigma,k_L}^*)^{-1} \\
&= s_{f^\sigma \circ f_\sigma} = s_{f_\sigma^+\circ f} \\
&= s_f \circ f_{k_L}^* \circ \alpha_{\cal{Y}^+}(\sigma) \circ (f_{k_L}^*)^{-1}
\end{align*}
as required.
\qed\medskip

In particular, under the assumption that $X$ has good reduction over $L$, the associated cohomology class
$$
\alpha^L_{X,\cal{L}} \,:=\, [\alpha_{\cal{Y}}] \in H^1(G,\cal{W}_{X_L,\cal{L}_L})
  $$
is independent of the chosen $\cal{L}_L$-terminal model $\cal{Y}$.

\begin{Proposition} 
 \label{prop: NOS1} 
 There exists an $\cal{L}_L$-terminal model of $X_L$ over $\OO_L$ for which the $G$-action is regular, if and only if the cohomology class
$$
\alpha^L_{X,\cal{L}} \,\in\, H^1\left(G,\cal{W}_{X_L,\cal{L}_L}\right)
 $$
is trivial. 
 In particular, $X$ has good reduction over $K$ if and only if $\alpha^L_{X,\cal{L}}$ is trivial.
\end{Proposition}

\prf 
For the first statement, one direction is clear, so suppose that the cohomology class $ \alpha^L_{X,\cal{L}}$ is trivial. 
Fix some $\cal{L}_L$-terminal model $\cal{Y}$, with associated cocycle $\alpha$. 
By assumption, there exists some $w\in \cal{W}_{X_L,\cal{L}_L}$ such that $\alpha(\sigma)= w^{-1}\sigma(w)$ 
for all $\sigma\in G$. 
We also know from Corollary \ref{cor: modweyl} that there exists some other terminal model 
$f:\cal{Y}\dashrightarrow \cal{Y}^+$ such that $s_f=w^{-1}$. 
Let 
$$
\alpha^+\,:\, G \,\to\, \cal{W}_{X_L,\cal{L}_L}
$$
be the cocycle associated to the model $\cal{Y}^+$. 
Then, using Proposition \ref{prop: cobound} we can see that $\alpha^+$ is trivial. 
Therefore, by Proposition \ref{prop: cocyc} we know that the rational $G$-action on $\cal{Y}^+$ is in fact regular.

For the second statement, the if direction follows from the first statement combined 
with Proposition \ref{prop: quotient}. 
Conversely, suppose that $X$ has good reduction over $K$, say with model $\cal{X}$. 
Then, by \cite[Proposition 4.5]{LM18} we may assume that $\cal{X}$ is $\cal{L}$-terminal, 
and hence $\cal{X}\otimes_{\OO_K} \OO_L$ is an $\cal{L}_L$-terminal model of $X_L$ 
to which the $G$-action extends. 
Thus, the cohomology class $\alpha^L_{X,\cal{L}}$ is trivial.
\qed\medskip

To give a condition that does not depend on the choice of unramified field extension $L/K$, we simply pass to the limit. 
So let 
$$ 
\cal{W}^\mathrm{nr}_{X,\cal{L}} \,:=\, 
  \mathrm{colim}_{K\subseteq L\subseteq K^\mathrm{nr}}\, \cal{W}_{X_L,\cal{L}_L}
  \,\cong\, \cal{W}_{X_{K^\mathrm{nr}},\cal{L}_{K^\mathrm{nr}}}
$$
be the colimit over all finite and unramified extensions of $K$ inside the fixed algebraic closure $\overline{K}$. Note that this direct system eventually stabilises (once all components of the exceptional divisor are defined over the residue field $k_L$ of $L$), so $\cal{W}^\mathrm{nr}_{X,\cal{L}}$ is the Weyl group of a Dynkin diagram of ADE-type. Then, there is a natural $G_k\cong\mathrm{Gal}(K^\mathrm{un}/K)$-action on $\cal{W}^\mathrm{nr}_{\cal{L}}$, 
and the above cohomology classes
$$
  \alpha^L_{X,\cal{L}} \,\in\, H^1\left(\mathrm{Gal}(L/K),\cal{W}_{X_L,\cal{L}_L}\right)
$$
for finite and unramified Galois extensions $L/K$ give rise to a well-defined class
$$
\alpha^\mathrm{nr}_{X,\cal{L}} \,\in\, H^1\left(G_k,\cal{W}^\mathrm{nr}_{X,\cal{L}}\right)
$$
in \emph{continuous} cohomology for the pro-finite group $G_k$. 
We can then rephrase Proposition \ref{prop: NOS1} as follows. 

\begin{Corollary} 
 \label{cor: goodred1} 
  Let $X/K$ be a K3 surface  satisfying $(\star)$ and the equivalent conditions of Theorem \ref{thm: first main result in introduction}. 
  Then, the following are equivalent.
 \begin{enumerate}  
   \item $X$ has good reduction over $K$.
   \item The cohomology class $\alpha^\mathrm{nr}_{X,\cal{L}}\in H^1(G_k,\cal{W}^\mathrm{nr}_{X,\cal{L}})$ is trivial for \emph{all} ample line bundles $\cal{L}$ on $X$.
   \item The cohomology class $\alpha^\mathrm{nr}_{X,\cal{L}}\in H^1(G_k,\cal{W}^\mathrm{nr}_{X,\cal{L}})$ is trivial for \emph{some} ample line bundle $\cal{L}$ on $X$.
 \end{enumerate} 
\end{Corollary}

\prf Since we have
$$
H^1\left(G_k,\cal{W}^\mathrm{nr}_{X,\cal{L}}\right) \,=\, 
\mathrm{colim}_{L}\, H^1\left(\mathrm{Gal}(L/K),\cal{W}_{X_L,\cal{L}_L}\right),
$$
where the colimit is taken over all finite and unramified Galois sub-extensions 
$K\subseteq L\subset K^\mathrm{nr}$, this follows from Proposition \ref{prop: NOS1}.
\qed\medskip

\section{Description of \texorpdfstring{$\Het{2}(X_{\overline{K}},\QQ_\ell)$}{H2et} and \texorpdfstring{$\DDcris(\Het{2}(X_{\overline{K}},\QQ_p))$}{Dcris}}

The purpose of this section is to prove Theorem \ref{theo: intro descr} from the introduction, which, for a K3 surface $X$ over $K$ satisfying $(\star)$ 
and the equivalent conditions of Theorem \ref{thm: first main result in introduction}, 
describes the unramified $G_K$-representation
$$ 
\Het{2}(X_{\overline{K}},\QQ_\ell),
$$
and, when $p>0$, the $F$-isocrystal
$$
\DDcris(\Het{2}(X_{\overline{K}},\QQ_p)),
$$
reasonably explicitly in terms of the $G_k$-representation $\Het{2}(Y_{\bar{k}},\QQ_\ell)$ and the $F$-isocrystal $\Hcris{2}(Y/K_0)$ 
associated to the canonical reduction $Y$ of $X$. 
To achieve this, we will first need a result on compatibility of comparison isomorphisms with cycle class maps, which briefly appeared earlier during the proof of Proposition \ref{prop: sf first props}. 

\subsection{Rational maps and crystalline comparison} 

Suppose that we have a birational map $f:\cal{X}\dashrightarrow \cal{X}^+$ between smooth models of K3 surfaces 
$X$ and $X^+$, respectively. 
Then, we may consider the graph of $f$
$$
\Gamma_f \,\subset\, \cal{X}\times_{\OO_K} \cal{X}^+
$$
as well as its generic and special fibers
$$
\Gamma_{f_K} \,\subset\, X\times_K X^+\mbox{ \quad and \quad } \Gamma_{f,k} \,\subset\, \cal{X}_k\times_k \cal{X}^+_k,
$$
respectively. 
Note that $\Gamma_{f_K}$ is genuinely the graph of the induced \emph{isomorphism} 
$f_K:X\rightarrow X^+$, but the same is not true of $ \Gamma_{f,k}$ in general. 
Since the cycle $\Gamma_{f,k}$ is purely of dimension $2$, we obtain maps on cohomology
$$
\begin{array}{ccccc}
   \Gamma_{f,k}^* &:& \Het{n}(\cal{X}^+_{\bar k},\QQ_\ell) &\rightarrow& \Het{n}(\cal{X}_{\bar k},\QQ_\ell) \\
     &&\alpha &\mapsto& p_{1*} \left( [\Gamma_{f,k}]\cup  p_2^*(\alpha)  \right) 
  \end{array}
$$
for $\ell\neq p$ and 
$$
\begin{array}{ccccc}
   \Gamma_{f,k}^*&:& \Hcris{n}(\cal{X}^+_k/K_0) &\rightarrow& \Hcris{n}(\cal{X}_k/K_0) \\
     && \alpha &\mapsto& p_{1*} \left(   [\Gamma_{f,k}]\cup  p_2^*(\alpha) \right)
\end{array}
$$
when $p>0$. 
The result we require is the following. 

\begin{Lemma}  
 \label{lemma: graphcomp} 
 In the above situation, the diagrams
$$ \xymatrix{   
    \Het{n}(X^+_{\overline{K}},\QQ_\ell)  \ar[r] \ar[d]_{f_K^*} & \Het{n}(\cal{X}^+_{\bar k} ,\QQ_\ell) \ar[d]^{\Gamma_{f,k}^*} \\  
    \Het{n}(X_{\overline{K}},\QQ_\ell ) \ar[r]  & \Het{n}(\cal{X}_{\bar k},\QQ_\ell)  }
 $$ 
for $\ell\neq p$, and 
$$\xymatrix{   
    \DDcris\left( \Het{n}(X^+_{\overline{K}},\QQ_p)\right) \ar[r] \ar[d]_{f_K^*} & \Hcris{n}(\cal{X}^+_k/K_0) \ar[d]^{\Gamma_{f,k}^*} \\  
    \DDcris\left( \Het{n}(X_{\overline{K}},\QQ_p)\right) \ar[r]  & \Hcris{n}(\cal{X}_k/K_0)  }
 $$
when $p>0$, commute. 
 Here, the horizontal arrows are the isomorphisms provided by the smooth and proper base change theorem and the 
 crystalline comparison theorem, respectively.
\end{Lemma}

\begin{Remark} 
 Since $\Het{n}(X_{\overline{\widehat{K}}},\QQ_p) \cong \Het{n}(X_{\overline{K}},\QQ_p)$ 
 as $G_{\widehat{K}}\cong G_K$ representations, there exists such a crystalline comparison theorem 
 without necessarily assuming $K$ to be complete.
\end{Remark}

\prf
Since the $\ell$-adic case for $\ell\neq p$ was handled in \cite[Lemma 5.6]{LM18}, we will only consider the 
$p$-adic case. As usual, we may assume that $K=\widehat{K}$ is complete. 
Since $\cal{X}$ and $\cal{X}^+$ have schematic fibers, the completions $\mathfrak{X}$ and $\mathfrak{X}^+$ along their special fibers are therefore smooth and proper formal schemes
(see the proof of Theorem \ref{thm:crysrep}).
Hence, we have Berthelot--Ogus comparison isomorphisms \cite[Theorem 2.4]{BO83}
$$
\begin{array}{lcl}
    \Hcris{n}(\cal{X}_k/W)\otimes_W K  &\cong& \HdR{n}(X/K)   \\ 
     \Hcris{n}(\cal{X}^+_k/W)\otimes_W K &\cong& \HdR{n}(X^+/K)
 \end{array}
$$
such that the diagram
$$ \xymatrix{   
  \DDcris\left( \Het{n}(X_{\overline{K}},\QQ_p)\right) \otimes_{K_0} K \ar[r]^-{\cong}\ar[d]_-{\cong}  
     & \Hcris{n}(\cal{X}_k/W)\otimes_W K \ar[d]^-{\cong}  \\ 
  \DDdR\left(\Het{n}(X_{\overline{K}},\QQ_p)\right)\ar[r]^-{\cong} 
     & H^n_\mathrm{dR}(X/K)      }
 $$
commutes (and similarly for $\cal{X}^+$), see for example \cite[Corollary 5.26]{CN17}. 
Since the induced map $f_K:X\dashrightarrow X^+$ is regular, it follows that the diagram
$$ \xymatrix{  
    \DDdR\left(\Het{n}(X_{\overline{K}},\QQ_p)\right)\ar[r]^-{\cong}\ar[d]_{\DDdR(f_K^*)} & \HdR{n}(X/K) \ar[d]^{f_K^*}  \\ 
    \DDdR\left(\Het{n}(X^+_{\overline{K}},\QQ_p)\right)\ar[r]^-{\cong} & \HdR{n}(X^+/K)   }
$$
commutes. 
It therefore suffices to prove that the Berthelot--Ogus comparison isomorphisms are compatible with 
$\Gamma_f^*$, in other words that the diagram
$$ 
   \xymatrix{  
     \Hcris{n}(\cal{X}_k/W)\otimes_W K \ar[r]^-{\cong}\ar[d]_{\Gamma_{f,k}^*} & \HdR{n}(X/K) \ar[d]^{f_K^*}  \\ 
     \Hcris{n}(\cal{X}^+_k/W)\otimes_W K\ar[r]^-{\cong} & \HdR{n}(X^+/K)   
    } 
$$
commutes. 
Since the horizontal isomorphisms are compatible with cup products (and hence Poincar\'e duality), 
it suffices to show that they are also compatible with cycle classes. 
In other words, we are given a smooth and proper algebraic space $\cal Z$ over $\OO_K$,
whose generic fiber $Z$ and whose special fiber $\cal{Z}_k$ are both schemes, 
and a closed and integral subspace $\cal T\subset \cal Z$ that is flat and of relative 
codimension $c$ over $\OO_K$, with generic fiber $T$ and special fiber $\cal{T}_k$. 
We must show that the isomorphism
$$
   H^{2c}_\mathrm{cris}(\cal{Z}_k/W)\otimes_W K \,\cong\, H^{2c}_\mathrm{dR}(Z/K) 
 $$
identifies $\mathrm{cl}(\cal{T}_k)\otimes 1$ with $\mathrm{cl}(T)$. Since $H^{2c}_\mathrm{cris}(\cal{Z}_k/W)\otimes_W K \cong H^{2c}_\mathrm{rig}(\cal{Z}_k/K)$ we can replace crystalline cohomology by rigid cohomology, and since rigid cohomology and algebraic de Rham cohomology both satisfy \'etale descent, by \cite{CT03} and \cite[Exp. VII, Proposition 4.3]{SGA4ii} respectively, the result follows by pulling back to an \'etale hypercover 
$\cal{Z}_\bullet\rightarrow \cal{Z}$ of $\cal{Z}$ by smooth $\OO_K$-schemes 
and applying \cite[Corollary 1.5.1]{CCM13}.
\qed\medskip

\subsection{Extending elements of \texorpdfstring{$\bm{\cal{W}}$}{W} to automorphisms of cohomology}

We can now use Lemma \ref{lemma: graphcomp} to give a different interpretation of the element $s_f\in \cal{W}_{X,\cal{L}}$ 
associated to a birational map $f:\cal{X}\dashrightarrow \cal{X}^+$ between terminal models of polarized K3 surfaces in 
Section \ref{sec: weyl k3}. 
This interpretation already appeared during the proof of Proposition \ref{prop: sf first props}, but it will be helpful to spell it out more 
clearly, including the $p$-adic case.

First of all, we will need to show how to extend elements of the Weyl group to automorphisms on cohomology. 
So, let $(X,\cal{L})$ be a polarized K3 surface over $K$, and $\cal{X}$ an $\cal{L}$-terminal model for $X$. 
We therefore have the Weyl group
$$
\cal{W}_{X,\cal{L}} \,\leq\, \mathrm{Aut}_{\ZZ}\left(\mathrm{Pic}(Y)\right)
$$
acting on the Picard group of the special fiber $Y$ of $\cal{X}$. 
Moreover, after base changing to $K^\mathrm{nr}$ we also have the corresponding 
$G_k$-equivariant version
$$
\cal{W}^\mathrm{nr}_{X,\cal{L}} \,\leq\, \mathrm{Aut}_{\ZZ}\left(\mathrm{Pic}(Y_{\bar k})\right)
$$
over $\bar k$.

\begin{Lemma}  
 \label{lemma: weyl to cohom} 
  There are $G_k$-equivariant and injective homomorphisms
  $$
   \begin{array}{lclll}
     i_\ell:\cal{W}_{X,\cal{L}}^\mathrm{nr} &\to& \mathrm{Aut}_{\QQ_\ell} & \left(\Het{2}(Y_{\bar{k}},\QQ_\ell(1)) \right)  & (\ell\neq p)\\
     i_p:\cal{W}_{X,\cal{L}}^\mathrm{nr} &\to& \mathrm{Aut}_{K_0^\mathrm{nr},F}&\left(\Hcris{2}(Y/K_0)(1) \otimes_{K_0} K_0^\mathrm{nr}\right) & (p>0).
   \end{array}
  $$
\end{Lemma}

\prf 
If $E\subset E_{X,\cal{L},\bar k}$ is an irreducible component of the geometric exceptional locus, then, 
for $\ell\neq p$, we simply take the reflection $s_E$ to the map
$$
\begin{array}{ccc}
   \Het{2}(Y_{\bar k},\QQ_\ell(1)) &\rightarrow& \Het{2}(Y_{\bar k},\QQ_\ell(1)) \\
   \alpha &\mapsto& \alpha + ([E]\cup \alpha )[E]
 \end{array}
$$
via the canonical identification $\Het{4}(Y_{\bar k},\QQ_\ell(2)) \cong \QQ_\ell$. 
When $p>0$, we choose a finite unramified extension $L/K$ with residue field $k_L$ over which $E$ is defined and 
send $s_E$ to the scalar extension of the map
$$
\begin{array}{ccc}
   \Hcris{2}(Y_{k_L}/L_0)(1) &\rightarrow& \Hcris{2}(Y_{k_L}/L_0)(1) \\
   \alpha &\mapsto& \alpha + ([E]\cup \alpha )[E]
  \end{array}
$$
where $L_0=W(k_L)[1/p]$, again via $\Hcris{4}(Y_{k_L}/L_0)(2) \cong L_0$.
\qed\medskip

In particular, by taking $G_k$-invariants we obtain injective homomorphisms
$$
\begin{array}{lclll}
  \cal{W}_{X,\cal{L}} &\to& \mathrm{Aut}_{G_k} & \left(\Het{2}(Y_{\bar{k}},\QQ_\ell(1)) \right)  & (\ell\neq p)\\
  \cal{W}_{X,\cal{L}} &\to& \mathrm{Aut}_{K_0,F}&\left(\Hcris{2}(Y/K_0)(1)\right) & (p>0).
\end{array}
$$
If we have a terminal birational map $f:\cal{X}\dashrightarrow \cal{X}^+$ to some other smooth model of a K3 surface $X^+/K$, then by \cite[Proposition 4.7]{LM18} $f$ is defined away from a finite collection of curves on the special fibers, and by minimality the induced rational maps
\begin{align*} f_K:X &\dashrightarrow X^+ \\
f_k: Y &\dashrightarrow \cal{X}^+_k
\end{align*}
on the generic and special fibers are isomorphisms. As well as the obvious pull-back maps
$$
\begin{array}{rlcll} 
	f_{k,\ell}^*:& \Het{2}(\cal{X}^+_{\bar k},\QQ_\ell(1)) & \to & \Het{2}(Y_{\bar k},\QQ_\ell(1)) & (\ell\neq p) \\
	f_{k,p}^* : & \Hcris{2}(\cal{X}^+_k/K_0)(1) & \to & \Hcris{2}(Y/K_0)(1) &(p>0)
\end{array}
$$
we also have ``generic fiber'' pull-back maps
$$
\begin{array}{rlcll} 
	f_{K,\ell}^*:& \Het{2}(\cal{X}^+_{\bar k},\QQ_\ell(1)) & \to & \Het{2}(Y_{\bar k},\QQ_\ell(1)) & (\ell\neq p) \\
	f_{K,p}^* : & \Hcris{2}(\cal{X}^+_k/K_0)(1) & \to & \Hcris{2}(Y/K_0)(1) &(p>0)
\end{array}
$$
which are defined by the commutative diagrams
$$\xymatrix{   
  \Het{2}(X^+_{\overline{K}},\QQ_\ell(1))  \ar[r]^-{\cong} \ar[d] & \Het{2}(\cal{X}^+_{\bar k} ,\QQ_\ell(1)) \ar[d]^{f_{K,\ell}^*} \\  
  \Het{2}(X_{\overline{K}},\QQ_p(1))  \ar[r]^-{\cong} & \Het{2}(Y_{\bar k},\QQ_\ell(1))} 
$$
for $\ell\neq p$, and
$$  \xymatrix{  
  \DDcris\left( \Het{2}(X^+_{\overline{K}},\QQ_p(1)) \right)  \ar[r]^-{\cong}\ar[d] & \Hcris{2}(\cal{X}^+_k/K_0)(1) \ar[d]^{f_{K,p}^*} \\  
  \DDcris\left( \Het{2}(X_{\overline{K}},\QQ_p(1)\right)  \ar[r]^-{\cong} & \Hcris{2}(Y/K_0)(1), } 
$$
when $p>0$. Here, the horizontal arrows are the appropriate comparison theorems, and the left hand vertical arrows simply pull-back on the generic fiber. 
We have the following corollary of Lemma \ref{lemma: graphcomp}.

\begin{Theorem} 
 \label{theo: sf diff} 
 Let $f:\cal{X}\dashrightarrow \cal{X}^+$ be a terminal birational map between terminal models of K3 surfaces over $K$. 
 Then, we have
$$
\begin{array}{llll} 
  i_\ell(s_f) = f^*_{K,\ell}\circ (f^*_{k,\ell})^{-1}:& \Het{2}\left(Y_{\bar k},\QQ_\ell(1)\right) &\to  \Het{2}\left(Y_{\bar k},\QQ_\ell(1)\right) & (\ell\neq p) \\
  i_p(s_f)   = f^*_{K,p}\circ (f^*_{k,p})^{-1} :& \Hcris{2}\left(Y/K_0\right)(1) &\to  \Hcris{2}\left(Y/K_0\right)(1) & (p>0). 
    \end{array}
 $$
\end{Theorem}

\prf
Let $\tilde{s}_f$ be the map $\mathrm{Pic}(\cal{X}^+_k)\rightarrow \mathrm{Pic}(Y)$ 
induced by the special fibre $\Gamma_{f,k}$ of the graph of $f$, 
and for all primes $\ell$, including $\ell=p$ when $p>0$, let $\tilde{s}_{f,\ell}$ be the map
$$
\begin{array}{rlcll} 
 \Gamma_{f,k}^*:& \Het{2}(\cal{X}^+_{\bar k},\QQ_\ell(1)) & \to & \Het{2}(Y_{\bar k},\QQ_\ell(1)) & (\ell\neq p) \\
 \Gamma_{f,k}^*: & \Hcris{2}(\cal{X}^+_k/K_0)(1) & \to & \Hcris{2}(Y/K_0)(1) &(p>0)
\end{array}
$$
induced in cohomology by the same cycle. 
Then, Lemma \ref{lemma: graphcomp} tells us that we have $f_{K,\ell}^*=\tilde{s}_{f,\ell}$, and since $s_f=\tilde{s}_f \circ (f_k^*)^{-1}$,  the result follows.
\qed\medskip

\subsection{Cohomological realisations of \texorpdfstring{$\bm{\alpha^L_{X,\cal{L}}}$}{al} and \texorpdfstring{$\bm{\alpha^\mathrm{nr}_{X,\cal{L}}}$}{anr}}

Now, let us suppose that $(X,\cal{L})$ is a polarized K3 surface over $K$, satisfying $(\star)$ and the equivalent conditions 
of Theorem \ref{thm: first main result in introduction}. 
Let $Y$ denote the canonical reduction of $X$. 

We will let $L/K$ denote a finite and unramified Galois extension over which $X$ has good reduction. 
Write $G$ for the Galois group of $L/K$ and $k_L/k$ for the residue field extension. 
We will write $L_0=W(k_L)[1/p]\subset \widehat{L}$ and $G_{k_L}$ for the absolute Galois group of $k_L$. We therefore have the exact sequence
$$
1 \,\to\,  G_{k_L} \,\to\, G_k \,\to\, G \,\to\, 1.
$$
Taking $G_{k_L}$-invariants in Lemma \ref{lemma: weyl to cohom} and twisting,  we obtain $G$-equivariant homomophisms
$$
\begin{array}{lclll}
  \cal{W}_{X_L,\cal{L}_L} &\to& \mathrm{Aut}_{G_{k_L}} & \left(\Het{2}(Y_{\bar{k}},\QQ_\ell(1)) \right)  & (\ell\neq p)\\
  \cal{W}_{X_L,\cal{L}_L} &\to& \mathrm{Aut}_{L_0,F}&\left(\Hcris{2}(Y_{k_L}/L_0)(1)\right) & (p>0).
\end{array}
$$

\begin{Definition} We define the \emph{$\ell$-adic and $p$-adic realisations} of $\alpha^L_{X,\cal{L}}$ to be
$$
\begin{array}{lllll}
 \beta^L_{X,\ell} &:= i_\ell(\alpha^L_{X,\cal{L}}) & \in& H^1\left(G,\mathrm{Aut}_{G_{k_L}}\left(\Het{2}(Y_{\bar k},\QQ_\ell)\right)\right) &(\ell\neq p) \\
 \beta^L_{X,p} &:= i_p(\alpha^L_{X,\cal{L}}) & \in & H^1\left(G,\mathrm{Aut}_{F,L_0}\left(\Hcris{2}(Y_{k_L}/L_0)\right)\right) &(p>0)
\end{array}
$$
 respectively.
\end{Definition}

We can now twist the cohomology groups $\Het{2}(Y_{\bar k},\QQ_\ell)$ and $\Hcris{2}(Y/K_0)$ 
by the classes $\beta^L_{X,\ell}$ and $\beta^L_{X,p}$, respectively, via the general procedure outlined 
in Section \ref{subsec: descent} and Section \ref{subsec: descent2}. 
We thus obtain new $G_k$-modules and $F$-isocrystals
$$
\Het{2}(Y_{\bar k},\QQ_\ell)^{\beta^L_{X,\ell}} \mbox{ \quad and \quad } \Hcris{2}(Y/K_0)^{\beta^L_{X,p}},
$$
respectively.

\begin{Theorem}  
 \label{theo: twist} 
 Let $(X,\cal{L})$ be a polarized K3 surface over $K$  satisfying $(\star)$ and the equivalent conditions of 
 Theorem \ref{thm: first main result in introduction}. 
 Then, there are natural isomorphisms 
 $$
  \begin{array}{lcll}
    \Het{2}(X_{\overline{K}},\QQ_\ell) &\cong& \Het{2}(Y_{\overline{k}},\QQ_\ell)^{\beta^L_{X,\ell}} & (\ell\neq p) \\
    \DDcris(\Het{2}(X_{\overline{K}},\QQ_p)) &\cong& \Hcris{2}(Y/K_0)^{\beta^L_{X,p}} & (p>0)
  \end{array}
 $$
of $G_k$-modules and $F$-isocrystals over $k$ respectively.
\end{Theorem}

\begin{Remark} This result implies that the cohomology classes $\beta_{X,\ell}^L$ and $\beta^L_{X,p}$ do not depend on the choice of polarisation $\mathcal{L}$.
\end{Remark}

\prf This is in fact reasonably straightforward. 
Suppose that we have a smooth, $\cal{L}_L$-terminal model $\cal{Y}$ for $X_L$ over $\OO_L$. This gives rise to a cospecialisation isomorphism
$$
\mathrm{comp}_{\cal{Y},\ell} \,:\, \Het{2}(X_{\overline{K}},\QQ_\ell) \,\isomto\,  \Het{2}(Y_{\overline{k}},\QQ_\ell),
$$
which is $G_{k_L}$-equivariant, although not $G_k$-equivariant in general. 
We can therefore define a function
$$
\beta_{\cal{Y},\ell} \,:\, G_k \rightarrow \mathrm{Aut}_{\QQ_\ell}\left(\Het{2}(Y_{\overline{k}},\QQ_\ell) \right)
$$
measuring the ``difference'' of the two actions. 
Concretely, we have
$$
\beta_{\cal{Y},\ell} (\sigma) \,=\, \mathrm{comp}_{\cal{Y},\ell}\circ  \sigma_{g}^*\circ \mathrm{comp}^{-1}_{\cal{Y},\ell} \circ  (\sigma_{s}^*)^{-1},
$$
where $\sigma_{g}^*$ is the action on $ \Het{2}(X_{\overline{K}},\QQ_\ell)$ and 
$\sigma_{s}^*$ that on $\Het{2}(Y_{\overline{k}},\QQ_\ell)$. 
Since $\mathrm{comp}_{\cal{Y},\ell}$ is $G_{k_L}$-equivariant, this map factors through $G$.
Moreover, by Theorem \ref{theo: sf diff} it coincides with the $\ell$-adic realization $i_{\ell} \circ \alpha_{\cal{Y}}$ 
of the cocycle
$$
\alpha_\cal{Y} :G \,\to\, \cal{W}_{X_L,\cal{L}_L}
$$
constructed in Section \ref{subsec: ogr}. 
In other words, $\mathrm{comp}_{\cal{Y}_\ell}$ becomes $G_k$-equivariant when we equip the left hand side with its 
natural action, and the right hand side with the ``twist'' of the natural action via the 
cocycle $i_{\ell} \circ \alpha_{\cal{Y}}=\beta_{\cal{Y},\ell}$, via the procedure described in Section \ref{subsec: descent2}.

In the $p$-adic case we argue similarly, using Theorem \ref{theo: sf diff} to show that the crystalline 
comparison isomorphism
$$
\left(\Het{2}(X_{\overline{L}},\QQ_p) \otimes_{\QQ_p} B_\mathrm{cris} \right)^{G_L} \,\isomto\, \Hcris{2}(Y_{k_L}/L_0)
$$
over $L$ becomes $G$-equivariant when we endow the left hand side with its natural $G$-action 
and the right hand side with the twist of the natural $G$-action via the cocycle 
$i_{p} \circ \alpha_{\cal{Y}}$. 
We can therefore apply Proposition \ref{prop: better descent}. 
\qed\medskip

We could equally well work over $\bar k$ as follows: using the $G_k$-equivariant homomorphisms
$$
\begin{array}{lclll}
  \cal{W}_{X,\cal{L}}^\mathrm{nr} &\to& \mathrm{Aut}_{\QQ_\ell} & \left(\Het{2}(Y_{\bar{k}},\QQ_\ell(1)) \right)  & (\ell\neq p)\\
  \cal{W}_{X,\cal{L}}^\mathrm{nr} &\to& \mathrm{Aut}_{K_0^\mathrm{nr},F}&\left(\Hcris{2}(Y/K_0)(1) \otimes_{K_0} K_0^\mathrm{nr}\right) & (p>0).
\end{array}
$$
 and taking a Tate twist, we obtain cohomology classes
$$
\begin{array}{ll}
 \beta^\mathrm{nr}_{X,\ell} := i_\ell(\alpha^\mathrm{nr}_{X,\cal{L}})  \in  H^1\left(G_k,\mathrm{Aut}_{\QQ_\ell}\left(\Het{2}(Y_{\bar k},\QQ_\ell)\right)\right) &(\ell\neq p) \\
 \beta^\mathrm{nr}_{X,p} := i_p(\alpha^\mathrm{nr}_{X,\cal{L}})  \in  H^1\left(G_k,\mathrm{Aut}_{F,K^\mathrm{nr}_0}\left(\Hcris{2}(Y/K_0)\otimes_{K_0}K_0^\mathrm{nr}\right)\right) &(p>0)
\end{array}
$$
that we can twist $\Het{2}(Y_{\bar k},\QQ_\ell)$ and $\Hcris{2}(Y/K_0)$ by, respectively.

\begin{Corollary}  
 \label{cor: twist} 
 Let $(X,\cal{L})$ be a polarized K3 surface over $K$  satisfying $(\star)$ and the equivalent conditions 
 of Theorem \ref{thm: first main result in introduction}. 
 Then, there are natural isomorphisms 
 $$
 \begin{array}{lcll}
    \Het{2}(X_{\overline{K}},\QQ_\ell) &\cong& \Het{2}(Y_{\overline{k}},\QQ_\ell)^{\beta^\mathrm{nr}_{X,\ell}} &(\ell\neq p) \\
    \DDcris(\Het{2}(X_{\overline{K}},\QQ_p)) &\cong& \Hcris{2}(Y/K_0)^{\beta^\mathrm{nr}_{X,p}} &(p>0)
  \end{array}
 $$
 of $G_k$-modules and $F$-isocrystals,respectively.
\end{Corollary}

\section{A N\'eron--Ogg--Shafarevich criterion for K3 surfaces}\label{sec: NOS}

We now come to the proof of our main result, namely Theorem \ref{theo: intro NOS} from the introduction. Of course, the implications (1)$\Rightarrow$(2) and (1)$\Rightarrow$(3) are just the smooth and 
proper base change and crystalline comparison theorems respectively, the real content is 
that (2)$\Rightarrow$(1) and (3)$\Rightarrow$(1). 
Given Corollary \ref{cor: goodred1} and Corollary \ref{cor: twist}, 
this amounts to showing that, given some ample line bundle $\cal{L}$ on $X$, the induced maps
$$
    i_\ell \,:\, H^1(G_k, \cal{W}^\mathrm{nr}_{X,\cal{L}}) \,\to \, H^1\left(G_k,\mathrm{Aut}_{\QQ_\ell}(\Het{2}(Y_{\bar k}, \QQ_\ell(1)))\right)
$$
for $\ell\neq p$, and 
 $$
  i_p\,: \, H^1(G_k, \cal{W}^\mathrm{nr}_{X,\cal{L}})  \,\to \, H^1\left(G_k,\mathrm{Aut}_{K^\mathrm{nr}_0,F}(\Hcris{2}(Y/K_0)(1)\otimes_{K_0} K^\mathrm{nr}_0)\right)
$$
when $p>0$, have trivial kernel. 
Equivalently, if we have some finite and unramified Galois extension $L/K$ with Galois group $G$ 
and residue field extension $k_L/k$, over which $X$ attains good reduction, it amounts to showing that the maps
$$
\begin{array}{lll}
  i_\ell : H^1(G, \cal{W}_{X_L,\cal{L}_L})  &\rightarrow H^1(G,\mathrm{Aut}_{G_{k_L}}(\Het{2}(Y_{\bar k}, \QQ_\ell(1)))) &(\ell\neq p)\\
  i_p: H^1(G, \cal{W}_{X_L,\cal{L}_L})  &\rightarrow H^1(G,\mathrm{Aut}_{L_0,F}(\Hcris{2}(Y_{k_L}/L_0)(1))) &(p>0)
\end{array}
$$
have trivial kernel. Note that we are free to enlarge $L$ if we wish, so we may assume that it is big enough so that $\mathcal{W}_{X_L,\mathcal{L}_L}=\mathcal{W}^\mathrm{nr}_{X,\mathcal{L}}$, that is, all geometric components of the exceptional divisor $E_{X,\mathcal{L}}$ are defined over $k_L$.
Let
$$
\begin{array}{lcll}  
 V_\ell &\subseteq& \Het{2}(Y_{\bar k}, \QQ_\ell(1)) & (\ell\neq p) \\
 V_p & \subseteq& \Hcris{2}(Y_{k_L}/L_0)(1)   & (p>0)
\end{array}
$$
denote the subspaces fixed by $G_{k_L}$ and Frobenius respectively.
These are therefore $G$-representations in finite dimensional $\QQ_\ell$ 
(resp. $\QQ_p$) vector spaces, and there are natural $G$-equivariant maps
$$
\begin{array}{llll}  \mathrm{Aut}_{G_{k_L}} & (\Het{2}(Y_{\bar k}, \QQ_\ell(1))) &  \rightarrow \mathrm{GL}(V_\ell) &(\ell\neq p) \\
 \mathrm{Aut}_{L_0,F} & (\Hcris{2}(Y_{k_L}/L_0)(1)) & \rightarrow \mathrm{GL}(V_p) & (p>0).
\end{array}
$$
It therefore suffices to show that the induced map
$$
i_\ell \,:\, H^1(G, \cal{W}_{X_L,\cal{L}_L})  \,\to\, H^1(G, \mathrm{GL}(V_\ell))
$$
has trivial kernel for all $\ell$ (including $\ell=p$ when $p>0$).

\begin{Proposition} 
  For all primes $\ell$ (including $\ell=p$ when $p>0$), if $\alpha\in H^1(G, \cal{W}_{X_L,\cal{L}_L})$ maps to the trivial class in
$$
H^1\left(G, \mathrm{GL}(V_\ell)\right), 
$$
then it maps to the trivial class in
$$
H^1\left(G, \mathrm{GL}(\Lambda_{X_L,\cal{L}_L,\QQ_\ell})\right).
 $$
\end{Proposition}

\begin{Remark} 
 Note that there is no naturally defined map between $ \mathrm{GL}(V_\ell)$ and $ \mathrm{GL}(\Lambda_{X_L,\cal{L}_L,\QQ_\ell})$
 in either direction, but the Weyl group $ \cal{W}_{X_L,\cal{L}_L}$ maps into both.
\end{Remark}

\prf  
This is entirely similar to the ``reduction to the transitive case'' part of the proof of Theorem \ref{theo: nagc main}. 
Consider the exact sequence
$$ 
0 \,\to\,  \Lambda_{X_L,\cal{L}_L,\QQ_\ell} \,\to\, V_\ell \,\to\,  T_\ell  \,\to\, 0
$$
of $G$-representations. This is also compatible with the action of $\cal{W}_{X_L,\cal{L}_L}$, and since all geometric components of the exceptional divisor $E_{X,\mathcal{L}}$ are defined over $k_L$, the induced action of the Weyl group 
on $T_\ell$ is trivial. Using the results of Section \ref{subsec: descent2}, we know that $\alpha$ being trivial in $H^1(G, \mathrm{GL}(V_\ell))$ 
is equivalent to the existence of an isomorphism $V_\ell^\alpha \cong V_\ell$ as $G$-representations, 
whereas $\alpha$ being trivial in $H^1(G, \mathrm{GL}( \Lambda_{X_L,\cal{L}_L,\QQ_\ell}))$ 
is equivalent to the existence of an isomorphism 
$\Lambda_{X_L,\cal{L}_L,\QQ_\ell}^\alpha \cong \Lambda_{X_L,\cal{L}_L,\QQ_\ell}$ as $G$-representations. 

Since the Weyl group acts trivially on $T_\ell$, we have $T_\ell^\alpha\cong T_\ell$, and therefore by semi-simplicity of the category of representations of the 
finite group $G$, we can conclude that $V_\ell^\alpha \cong V_\ell$ if and only if 
$\Lambda_{X_L,\cal{L}_L,\QQ_\ell}^\alpha \cong \Lambda_{X_L,\cal{L}_L,\QQ_\ell}$.
\qed\medskip

Thus, Theorem \ref{theo: intro NOS} amounts to showing that
$$
H^1\left(G, \cal{W}_{X_L,\cal{L}_L}\right) \,\to\, H^1(G, \mathrm{GL}\left(\Lambda_{X_L,\cal{L}_L,\QQ_\ell})\right)
$$
has trivial kernel. 
Since $G$ is finite, this is exactly the content of Theorem \ref{theo: nagc main}.

\section{Integral \texorpdfstring{$p$}{p}-adic Hodge theory}

So far we have worked with the rational \'etale cohomology of a K3 surface $X$ over $K$
with potentially good reduction.
In this last section, we will explain how all the results carry over to integral \'etale cohomology, 
both $\ell$-adically and $p$-adically. 

Assume therefore that $X$ satisfies $(\star)$ and the equivalent conditions of 
Theorem \ref{thm: first main result in introduction}. 
Then, we have maps 
$$
\begin{array}{lclll}
  \cal{W}_{X,\cal{L}}^\mathrm{nr} &\to& \mathrm{Aut}_{\QQ_\ell} & \left(\Het{2}(Y_{\bar{k}},\QQ_\ell) \right)  & (\ell\neq p)\\
  \cal{W}_{X,\cal{L}}^\mathrm{nr} &\to& \mathrm{Aut}_{K_0^\mathrm{nr},F}&\left(\Hcris{2}(Y/K_0) \otimes_{K_0} K_0^\mathrm{nr}\right) & (p>0).
\end{array}
$$
and it follows from their definition (Lemma \ref{lemma: weyl to cohom}) that they factor through maps
$$
\begin{array}{lclll}
   i_\ell:\cal{W}_{X,\cal{L}}^\mathrm{nr} &\to& \mathrm{Aut}_{\ZZ_\ell} & \left(\Het{2}(Y_{\bar{k}},\ZZ_\ell) \right)  & (\ell\neq p)\\
  i_p:\cal{W}_{X,\cal{L}}^\mathrm{nr} &\to& \mathrm{Aut}_{W^\mathrm{nr},F}&\left(\Hcris{2}(Y/W) \otimes_{W} W^\mathrm{nr}\right) & (p>0),
\end{array}
$$
where $W^\mathrm{nr}$ denote the ring of integers of $K_0^\mathrm{nr}$. 
Thus, we could equally well define cohomology classes
$$
\begin{array}{lll}
 \beta^\mathrm{nr}_{X,\ell} := i_\ell(\alpha^\mathrm{nr}_{X,\cal{L}}) & \in H^1(G_k,\mathrm{Aut}_{\ZZ_\ell}(\Het{2}(Y_{\bar k},\ZZ_\ell))) &(\ell\neq p) \\
 \beta^\mathrm{nr}_{X,p} := i_p(\alpha^\mathrm{nr}_{X,\cal{L}}) & \in H^1(G_k,\mathrm{Aut}_{F,W^\mathrm{nr}}(\Hcris{2}(Y/W)\otimes_{W}W^\mathrm{nr})) &(p>0)
\end{array}
$$
on the integral level, and we can twist
$$
\Het{2}(Y_{\bar k},\ZZ_\ell)\mbox{ \quad and \quad }\Hcris{2}(Y/W)
$$
by these classes to obtain new $G_k$-modules and $F$-crystals
$$
\Het{2}(Y_{\bar k},\ZZ_\ell)^{\beta^\mathrm{nr}_{X,\ell}} \mbox{ \quad and \quad }\Hcris{2}(Y/W)^{\beta^\mathrm{nr}_{X,p}},
$$
respectively. 
\medskip

\prf[Proof of Theorem \ref{thm: last main result in introduction}] 
The second statement is a trivial consequence of Theorem \ref{theo: intro NOS} and the first follows easily from the proof of 
Theorem \ref{theo: twist} as follows:
we need to show that if we choose an unramified extension and Galois extension $L/K$ with Galois group $G$ and residue field $k_L$, large enough so that $X$ has good reduction over $L$, 
then the induced comparison isomorphisms
$$
\begin{array}{lcll} 
\Het{2}(X_{\overline{K}},\ZZ_\ell) &\cong  &\Het{2}(Y_{\bar{k}},\ZZ_\ell) & (\ell\neq p) \\
\mathrm{BK}_{\OO_L} \left( \Het{2}(X_{\overline{L}},\ZZ_p) \right) \otimes_{\mathfrak{S},\sigma} W(k_L) &\cong & \Hcris{2}(Y_{k_L}/W(k_L))  & (p>0)
\end{array}
$$
provided by the smooth and proper base change theorem and Theorem \ref{thm:crysrep},
respectively, identify the natural $G_k$-action (resp. $G$-action) on the left hand side with the 
``$\beta$-twisted'' $G_k$-action (resp. $G$-action) on the right hand side. 
But we can embed the left hand side equivariantly in 
$\Het{2}(X_{\overline{K}},\QQ_\ell)$ (resp. $(\Het{2}(X_{\overline{L}},\QQ_p)\otimes_{\QQ_p} B_\mathrm{cris})^{G_L}$), 
and the right hand side equivariantly in $\Het{2}(Y_{\bar{k}},\QQ_\ell)$ (resp. $\Hcris{2}(Y_{k_L}/L_0)$). The result therefore follows from the proof of Theorem \ref{theo: twist}. 
\qed\medskip

\appendix

\addtocontents{toc}{\setcounter{tocdepth}{0}}

\section{Calculating the kernel of \texorpdfstring{$H^1(\ZZ/3,\cal{W}_{D_4}) \rightarrow H^1(\ZZ/3,\mathrm{GL}_4(F))$}{H1} }  \label{app: code D4}

\addtocontents{toc}{\setcounter{tocdepth}{1}}

Here, we reproduce some \cite{SAGE} code that calculates an upper bound for the size of the 
kernel of the natural map
$$
H^1\left(\ZZ/3,\cal{W}_{D_4}\right) \,\to\, H^1\left(\ZZ/3,\mathrm{GL}_4(F)\right),
$$
for the $\ZZ/\ZZ3$-action on $D_4$ and $F$ a field of characteristic $0$. 
It is explained during the proof of Proposition \ref{prop: nagc faith}, and its output is also included below. A SAGE file can be downloaded from \url{https://staff.fnwi.uva.nl/c.d.lazda/pubs/}.

\begin{lstlisting}[language=Python]
# define the basic reflections s1,s2,s3,s4
s=[];
s.append(matrix([[-1,1,0,0],[0,1,0,0],
		[0,0,1,0],[0,0,0,1]]));
s.append(matrix([[1,0,0,0],[1,-1,1,1],
		[0,0,1,0],[0,0,0,1]]));
s.append(matrix([[1,0,0,0],[0,1,0,0],
		[0,1,-1,0],[0,0,0,1]]));
s.append(matrix([[1,0,0,0],[0,1,0,0],
		[0,0,1,0],[0,1,0,-1]]));

# create the Weyl group as a set of matrices
# 192 is the size of the Weyl group in this case
W = [matrix.identity(4)];
while len(W) < 192:		
    for i in range(4):
        for w in W:
            if w*s[i] not in W:
                W.append(w*s[i]);
            if len(W) == 192:
                break;
        if len(W) == 192:
            break;

# check the Weyl group is actually a group
for i in range(4):
    for w in W:
        assert w*s[i] in W;
for w in W:
    assert w.inverse() in W;
assert matrix.identity(4) in W;

# define the matrix giving the action of a generator
# of Z/3 on the root lattice
g = matrix([[0,0,0,1],[0,1,0,0],[1,0,0,0],[0,0,1,0]]);

# check that 1 is an eigenvalue of g
assert 1 in g.eigenvalues(); 

# calculate the set of cocycles
cocycle = [];
for w in W:
    if w*g*w*g*w*g == matrix.identity(4):
        cocycle.append(w);

# calculate a representing set of cohomology classes
cohom = [];
for w in cocycle:
    i=0
    for r in W:
        i=i+1;
        M = r.inverse()*w*g*r*g.inverse();
        if M in cohom:
            break;
    if i == len(W):
        cohom.append(w);

# calculate an upper bound for the kernel by
# throwing out all w such that w*g does not
# have 1 as an eigenvalue
kernel=[];
for w in cohom:
    N = w*g;
    if 1 in N.eigenvalues():
        kernel.append(w);    

print('# kernel <= %d' % len(kernel));        
        
# prints the non-trivial cohomology classes        
print('non trivial cohomology class representatives:')
for w in cohom:
    if w != matrix.identity(4):
        w;
\end{lstlisting}
Output:
\begin{lstlisting}
# kernel <= 1
non trivial cohomology class representatives:
[-1  0  0  0]
[-1  0  1 -1]
[-1  1  0 -1]
[-1  0  1  0]
\end{lstlisting}

\bibliographystyle{../../Templates/bibsty}
\bibliography{../../lib.bib}

\end{document}